\documentclass[11pt]{article}
\usepackage[english]{babel} 
\usepackage{lineno}
\usepackage{graphicx,multicol}
\usepackage{epic,eepic,epsfig}
\usepackage{caption}

\usepackage{amsmath,amsfonts,amsthm,amssymb}
\usepackage{pstricks,pstricks-add,pst-math,pst-xkey}

\usepackage{algorithm, algorithmic}

\usepackage{graphicx,pstricks,pst-node,amssymb}

\usepackage{tikz, enumitem, caption, subcaption, float}
\usepackage{indentfirst, hyperref}
\usepackage{amsmath,amsfonts,amssymb,amsthm,epsfig,epstopdf,titling,url,array, mathdots}                     
\usepackage[margin=1in]{geometry}
\usepackage{mathtools}
                        
\usepackage{nccmath}    
\newtagform{roman}[]()
                        
\usetikzlibrary{snakes} 
\usetikzlibrary{decorations.markings}
\tikzset{unlabelled/.style={black, draw=black, circle, fill, scale=0.5}}
\tikzset{labelled/.style={black, draw=black, circle, scale=0.7, font=\Large}}
         
\tikzset{->-/.style={decoration={
                        markings,
                        mark=at position #1 with {\arrow{>}}},postaction={decorate}}}
\tikzset{dashedarc/.style={decoration={
                        markings,
                        mark=at position #1 with {\arrow{>}}},postaction={decorate}}}
                        
\usepackage{tikz}
\pgfdeclarelayer{edgelayer}
\pgfdeclarelayer{nodelayer}
\pgfsetlayers{edgelayer,nodelayer,main}

\tikzstyle{none}=[inner sep=0pt]
\definecolor{hexcolor0xf81e1c}{rgb}{0.973,0.118,0.110}
\definecolor{hexcolor0x3c00ff}{rgb}{0.235,0.000,1.000}

\tikzstyle{whitevertex}=[circle,fill=white,draw=black, scale = 0.5]
\tikzstyle{redvertex}=[circle,fill=hexcolor0xf81e1c,draw=black, scale = 0.5]
\tikzstyle{bluevertex}=[circle,fill=hexcolor0x3c00ff,draw=black, scale = 0.5]
\tikzstyle{greenvertex}=[circle,fill=green,draw=black, scale=0.5]
\tikzstyle{purplevertex}=[circle,fill=magenta,draw=black, scale=0.5]
\tikzstyle{grayvertex}=[circle,fill=white,draw=gray, scale=0.5]
\tikzstyle{blackvertex}=[circle,fill=black,draw=black, scale=0.5]

\tikzstyle{textbox}=[rectangle,fill=none,draw=none]
\tikzstyle{box}=[rectangle,fill=none,draw=black]

\tikzstyle{arc}=[black, ->]
\tikzstyle{grayarc}=[gray, ->]
\tikzstyle{bluearc}=[blue, ->]
\tikzstyle{grayedge}=[draw=gray]
\tikzstyle{blueedge}=[draw=blue]
\tikzstyle{rededge}=[draw=red]
\tikzstyle{edge}=[draw=black]

\tikzstyle{vertex}=[circle, ,fill=white,draw=black, scale=0.5]
\tikzstyle{10circle}=[circle, scale=10.0,draw=black]
\tikzstyle{10oval}=[ellipse, scale=10.0,draw=black]

\usepackage{amsmath,amsfonts,amsthm,amssymb}
\usepackage{pstricks,pstricks-add,pst-math,pst-xkey}

\usepackage{algorithm, algorithmic}

\usepackage{graphicx,pstricks,pst-node,amssymb}

\begin{document}

\newtheorem{tm}{\hspace{5mm}Theorem}
\newtheorem{claim}{\hspace{5mm}Claim}
\newtheorem{prp}[tm]{\hspace{5mm}Proposition}
\newtheorem{dfn}[tm]{\hspace{5mm}Definition}
\newtheorem{lemma}[tm]{\hspace{5mm}Lemma}
\newtheorem{cor}[tm]{\hspace{5mm}Corollary}
\newtheorem{conj}[tm]{\hspace{5mm}Conjecture}
\newtheorem{prob}[tm]{\hspace{5mm}Problem}
\newtheorem{remark}[tm]{\hspace{5mm}Remark}
\newtheorem{quest}[tm]{\hspace{5mm}Question}
\newtheorem{alg}[tm]{\hspace{5mm}Algorithm}
\newtheorem{sub}[tm]{\hspace{5mm}Algorithm}
\newcommand{\la}{\langle}
\newcommand{\ra}{\rangle}
\newcommand{\pf}{{\bf Proof: }}
\newcommand{\ppf}{{\bf Proof:}}
\newcommand{\beq}{\begin{equation}}
\newcommand{\eeq}{\end{equation}}
\newcommand{\<}[1]{\left\langle{#1}\right\rangle}

\newcommand{\Z}{\mathbb{$Z$}}
\newcommand{\Q}{\mathbb{$Q$}}
\newcommand{\R}{\mathbb{$R$}}

%\renewcommand{\thefootnote}{\fnsymbol{footnote}}

%%%%%%%%%%%%%%%%%%%%JBJ%%%%%%%%%%%%%%%%%%%%%%%%%%%%%%%%%%
\newcommand{\jbj}[1]{{\color{blue}#1}}
\newcommand{\obs}[1]{{\color{magenta}#1}}
%\usepackage{eepic,epsfig}

%\theoremstyle{definition}
%\newtheorem{remark}[tm]{Remark}
%\newtheorem{exercise}{Exercise}[section]
%\newtheorem{exercice}[exercise]{Exercise}
%\newtheorem{algo}{Algorithm}[section]
%\newtheorem{problem}[tm]{Problem}
%\newtheorem{conjecture}[tm]{Conjecture}
%\newtheorem{definition}[tm]{Definition}
%\newtheorem{question}[tm]{Question}

%---------------------------------------------------------------------

%--------- The following lines are needed to fill out the page

%---------------------------------------------------------------------

\setlength{\topmargin}{0cm}

\setlength{\headheight}{0.1cm}

\setlength{\headsep}{0.1cm}

\setlength{\textheight}{23.4cm}

\setlength{\oddsidemargin}{0.1cm}

\setlength{\evensidemargin}{0.1cm}

\setlength{\textwidth}{15.5cm}

% ---------------------------------------------------------------------

\newcommand{\perhapsREPLACE}[1]{{\color{blue} #1}}

\newcommand{\AlternateProof}[1]{{\color{red} #1}}

\newcommand{\dXY}[2]{d(#1,#2)}  % The number of arcs from #1 to #2

\newcommand{\ssX}{\hspace{1.2cm}}

\newcommand{\todo}[1]{\marginpar{{\textbf{JBJ}}\\\raggedright #1}}

\newcommand{\todoF}[1]{\marginpar{{\textbf{Fred}}\\\raggedright {\color{red} #1}}}
\newcommand{\LP}[1]{{\color{red} #1}}

\newcommand{\induce}[2]{\mbox{$ #1 \langle #2 \rangle$}}

\newcommand{\nsdom}{\mbox{$\not\Rightarrow$}}

\newcommand{\dom}{\mbox{$\rightarrow$}}

\newcommand{\sdom}{\mbox{$\mapsto$}}
\newcommand{\nad}{\ \mbox{${ \not\sim }$}\ }

\newcommand{\2}{\vspace{2mm}}

\theoremstyle{plain}

\newenvironment{subproof}{\par\noindent {\it Proof}.\ }{\hfill$\lozenge$\par\vspace{11pt}}

\bibliographystyle{plain}
%\pagewiselinenumbering
%\setpagewiselinenumbers
%\modulolinenumbers[1]
%\linenumbers

\title{On graphs which have locally complete 2-edge-colourings and 
         their relationship to proper circular-arc graphs}
\author{J\o{}rgen Bang-Jensen\thanks{Department of Mathematics and Computer
    Science, University of Southern Denmark, Odense, Denmark (email:
    jbj@imada.sdu.dk). Research supported by the Independent
    Research Fund Denmark under grant number DFF 7014-00037B}\and Jing Huang
\thanks{Department of Mathematics and Statistics,
University of Victoria, Victoria, B.C., Canada (email: huangj@uvic.ca). Research
supported by NSERC}}
\date{}

\maketitle

\begin{abstract}
A 2-edge-coloured graph $G$ is called {\bf locally complete} if for each vertex $v$,
the vertices adjacent to $v$ through edges of the same colour induce a complete
subgraph in $G$. Locally complete 2-edge-coloured graphs have nice properties and 
there exists a polynomial algorithm to decide whether such a graph has an alternating
hamiltonian cycle, where alternating means that the colour of two consecutive edges 
on the cycle are different. In this paper we show that graphs having locally 
complete 2-edge-colourings can be recognized in polynomial time. We give 
a forbidden substructure characterization for this class of graphs analogous 
to Gallai's characterization for cocomparability graphs. Finally, we characterize 
proper interval graphs and proper circular-arc graphs which have locally complete
2-edge-colourings by forbidden subgraphs.\\
\noindent{}{\bf Keywords:} Locally complete 2-edge-coloured graph, forbidden 
structure, characterization, recognition, proper interval graph,
proper circular-arc graph, local tournament
\end{abstract}

\section{Introduction}

%\obs{Do we want to use the name twins or true twins?}

The structure of a graph is often inherent in its ability of admitting certain
orientations or vertex orderings which are particularly useful in the design of 
efficient algorithms for solving hard problems in graphs. For example, comparability
graphs are precisely those which admit transitive orientations or quasi-transitive
orientations \cite{bh,gallai,ghouri,hh}, and chordal graphs are just those which
admit perfect elimination orderings \cite{fg,golumbic}. Other graph and digraph
classes with similar structure properties can be found in
\cite{dama,farber,golumbic,gg,hhl,hhmr,hmr,huangGallai,huang,huang-18,skrien}. 
In this paper we study graphs which admit certain edge-colourings.

All graphs considered in this paper are simple, i.e., they contain no loops or
multiple edges. We shall consider graphs whose edges are coloured with
two colours. These colourings are not necessarily proper in the sense that two edges
which share a common vertex may be coloured by the same colour. We refer to such
a graph as a {\bf 2-edge-coloured} graph. Many problems on paths and cycles in directed graphs are special cases of analogous problems on cycles in 2-edge-coloured graphs with the property that the edges on the cycle or path alternate between the two colours. There is a rich literature on alternating paths and cycles in 2-edge-coloured graphs, see e.g \cite{bang2009,bangDM165}.

A 2-edge-coloured graph $G$ is called {\bf locally complete} if for each vertex $v$,
the vertices adjacent to $v$ through edges of the same colour induce a complete
subgraph in $G$. Equivalently, $G$ contains no monochromatic induced path of length
2.

Locally complete 2-edge-coloured graphs were introduced in \cite{galeana-sanchez}
under the name of {\bf 2-M-closed} graphs. The authors of \cite{galeana-sanchez}
characterized locally complete 2-edge-coloured graphs which have alternating
Hamiltonian cycles. See also \cite{bangGC37} for further results on locally complete
2-edge-coloured graphs.

The definition of locally complete 2-edge-coloured graphs echoes another notion on
digraphs. A digraph $D$ is called {\bf locally semicomplete} if for each vertex $v$,
the vertices connecting to $v$ by out-going arcs from $v$ induce a semicomplete
subdigraph and the vertices connecting to $v$ by in-coming arcs to $v$ also
induce a semicomplete subdigraph in $D$. Local semicomplete digraphs are a natural
generalization of semicomplete digraphs and tournaments \cite{jbj}. When a local
semicomplete digraph is an oriented graph, it is called a {\bf local tournament}
\cite{huang}.

The underlying graphs of local tournaments are intimately related to proper
circular-arc graphs.
A graph is called a {\bf proper circular-arc} graph if it is the intersection graph
of a family of circular-arcs on a circle where no arc is properly contained in
another arc. Skrien \cite{skrien} proved that a connected graph is a proper
circular-arc graph if and only if it can be oriented as a local tournament.
Tucker \cite{tucker} gave a forbidden subgraph characterization of proper
circular-arc graphs. A graph is called a {\bf proper interval graph} if it is 
the intersection of a family of intervals in a line where no interval is properly 
contained in another interval. Proper interval graphs are a subclass of proper 
circular-arc graphs. In fact, they are precisely the graphs which can be oriented 
as non-strong local tournaments, cf. \cite{huang}. 
A description of all possible local tournament orientations of a proper 
circular-arc graph and all possible non-strong local tournament orientations of 
proper intervals graphs has been given in \cite{huang}.

In this paper, we study the graphs which have locally complete 2-edge-colourings.
We call such graphs {\bf locally complete 2-edge-colourable}. The complement 
$\overline{H}$ of any bipartite graph $H$ is locally complete 2-edge-colourable. 
Indeed, if $(X,Y)$ is a bipartition of $H$, then colouring all edges of 
$\overline{H}$ with both end vertices in $X$ or in $Y$ by one colour and the 
remaining edges with the other colour gives a locally complete 2-edge-colouring
of $\overline{H}$.

Clearly, if a graph $G$ is locally complete 2-edge-colourable then every induced
subgraph of $G$ is also locally complete 2-edge-colourable. Since a locally 
complete 2-edge-colouring of a graph does not contain a monochromatic induced 
path of length 2, the claw $K_{1,3}$ is not locally complete 2-edge-colourable.
In fact, no odd cycle of length $\geq 5$ or its complement is locally 
complete 2-edge-colourable. Hence locally complete 2-edge-colourable graphs are 
claw-free and perfect.

We show that a graph $G$ is locally complete 2-edge-colourable  if and only if its
associated auxiliary graph (see definition in Section~\ref{sec2}) is bipartite. This
leads to a polynomial time recognition algorithm and a forbidden structure
characterization for locally complete 2-edge-colourable graphs.

We compare locally complete 2-edge-colourable graphs with the graphs which can be 
oriented as local tournaments, and in particular, with proper circular-arc graphs. 
A complete list of minimal locally complete 2-edge-colourable graphs which are
not proper circular-arc graphs can be derived from Tucker's characterization of 
proper circular-arc graphs.
We also find the complete list of minimal proper circular-arc graphs which are not
locally complete 2-edge-colourable. Combining this with Tucker's forbidden subgraph 
characterization of proper circular-arc graphs we obtain a forbidden subgraph 
characterization of proper circular-arc graphs which are locally complete 
2-edge-colourable.

The paper is organized as follows. In Section \ref{sec2} we show that graphs which have a locally complete 2-edge-colouring can be recognized in polynomial time. We also give a forbidden structure characterization of graphs which have a locally complete 2-edge-colouring. In Section \ref{sec3} we characterize those graphs which have a locally complete 2-edge-colouring but are not proper circular-arc graphs. In Section \ref{sec4} we first characterize proper interval graphs which have a locally complete 2-edge-colouring and use that result to  characterize those graphs which have a locally complete 2-edge-colouring
 which are proper circular-arc graphs but are not proper interval graphs.

\section{Characterizations and recognition} \label{sec2}

Let $G$ be a graph. Define the {\bf auxiliary} graph $G^+$ of $G$ as follows:
the vertices of $G^+$ consists of all edges of $G$ and two vertices of $G^+$ are
adjacent if and only if the two corresponding edges form an induced path in $G$.
Any edge-colouring of $G$ transforms directly to a vertex-colouring of $G^+$ and
vice versa. If a 2-edge-colouring of $G$ is locally complete then the transformed
vertex-colouring of $G^+$ is a proper 2-colouring.
Conversely, if a 2-vertex-colouring of $G^+$ is proper then the transformed
edge-colouring of $G$ is locally complete. Hence we have the following:

\begin{tm} \label{aux}
A graph $G$ is locally complete 2-edge-colourable if and only if the auxiliary graph
$G^+$ is bipartite.
\qed
\end{tm}

It follows from Theorem \ref{aux} that to determine whether a graph $G$ has a locally
complete 2-edge-colouring one simply constructs the auxiliary graph $G^+$ and checks
whether it is bipartite or has a 2-vertex-colouring. Moreover, any 2-vertex-colouring
of $G^+$ (if exists) transforms to a locally complete 2-edge-colouring of $G$. 

\begin{cor}
There is a polynomial-time algorithm to decide if a graph has a locally complete
2-edge-colouring and find one if it exists. 
\qed
\end{cor}

\begin{figure}[ht]
\begin{center}
\begin{tikzpicture}[>=latex]
%\begin{pgfonlayer}{nodelayer}
\node [label={above:$v_0$}] [style=blackvertex] (0) at (0,1.5) {};
\node [label={right:$v_1$}] [style=blackvertex] (1) at (1.1,.7) {};
\node [label={right:$v_2$}] [style=blackvertex] (2) at (1.1,-.7) {};
\node [label={below:$v_3$}] [style=blackvertex] (3) at (0,-1.5) {};
\node [label={left:$v_4$}] [style=blackvertex] (4) at (-1.1,-.7) {};
\node [label={left:$v_5$}] [style=blackvertex] (5) at (-1.1,.7) {};

 \draw[dotted] (0) to (2) to (4) to (0);
 \draw[dotted] (1) to (3) to (5) to (1);

%\node [style=textbox]  at (0, .5) {\tiny{$W_5$}};

 \draw[dashed] (0) to (1) to (2) to (3) to (4) to (5) to (0);

%\node [style=circle, scale=.1] (01) at (.2,.2) {};
 \draw[dashed] (0) to (.2,.9);
 \draw[dashed] (0) to (-.2,.9);
 \draw[dashed] (0) to (0,.9);
 \draw[dashed] (1) to (.7,.4);
 \draw[dashed] (1) to (.7,.6);
 \draw[dashed] (2) to (.7,-.4);
 \draw[dashed] (2) to (.7,-.6);
 \draw[dashed] (3) to (.2,-.9);
 \draw[dashed] (3) to (-.2,-.9);
 \draw[dashed] (3) to (0,-.9);
 \draw[dashed] (4) to (-.7,-.4);
 \draw[dashed] (4) to (-.7,-.6);
 \draw[dashed] (5) to (-.7,.4);
 \draw[dashed] (5) to (-.7,.6);

\node [label={above:$v_0$}] [style=blackvertex] (6) at (6,2) {};
\node [label={right:$v_1$}] [style=blackvertex] (7) at (7.4,.9) {};
\node [label={right:$v_2$}] [style=blackvertex] (8) at (7.7,-.7) {};
\node [label={below:$v_3$}] [style=blackvertex] (9) at (6.8,-2) {};
\node [label={below:$v_4$}] [style=blackvertex] (10) at (5.2,-2) {};
\node [label={left:$v_5$}] [style=blackvertex] (11) at (4.3,-.7) {};
\node [label={left:$v_6$}] [style=blackvertex] (12) at (4.6,.9) {};

 \draw[dotted] (6) to (8) to (10) to (12) to (7) to (9) to (11) to (6);

 \draw[dashed] (6) to (7) to (8) to (9) to (10) to (11) to (12) to (6);

 \draw[dashed] (6) to (6.2,1.2);
 \draw[dashed] (6) to (5.8,1.2);
 \draw[dashed] (6) to (6,1.2);

 \draw[dashed] (7) to (7,.8);
 \draw[dashed] (7) to (7.1,.6);

 \draw[dashed] (8) to (7.3,-.5);
 \draw[dashed] (8) to (7.2,-.7);

 \draw[dashed] (9) to (6.5,-1.5);
 \draw[dashed] (9) to (6.8,-1.4);

 \draw[dashed] (10) to (5.2,-1.4);
 \draw[dashed] (10) to (5.5,-1.5);

 \draw[dashed] (11) to (4.7,-.5);
 \draw[dashed] (11) to (4.7,-.7);

 \draw[dashed] (12) to (5,.8);
 \draw[dashed] (12) to (4.9,.5);

%\node [style=textbox]  at (-4, -.6) {$(i)$};
%\node [style=textbox]  at (-.5, -.6) {$(ii)$};
%\node [style=textbox]  at (3.5, -.6) {$(iii)$};
%\node [style=textbox]  at (7.5, -.6) {$(iv)$};

%\end{pgfonlayer}
\end{tikzpicture}
\end{center}
\vspace{-2mm}
\caption{Kaleidoscopes of order 6 (left) and 7 (right) where dotted lines are
walks and dashed lines are non-edges.}
\label{kaleidoscopes}
\end{figure}
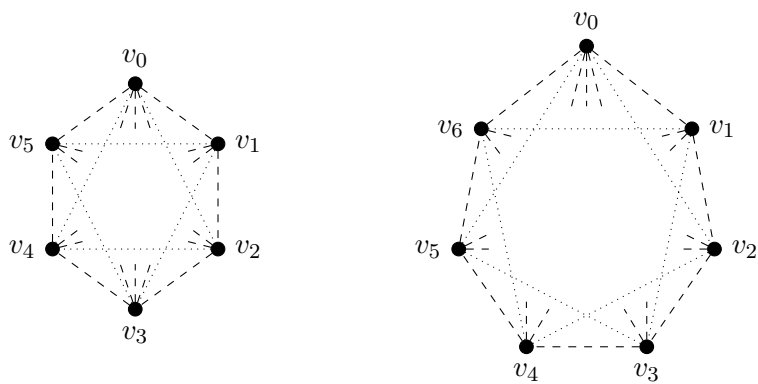

Theorem \ref{aux} can be used to describe a forbidden structure for locally complete 2-edge-colourable graphs.
Let $H$ be a graph. A walk $W$ {\bf avoids} a vertex $v$ in $H$ if $W$ does not 
contain $v$ or any neighbour of $v$. A {\bf kaleidoscope} of order $k$ ($k \geq 2$) 
in $H$ consists of (not necessarily distinct) vertices $v_0, v_1, \dots, v_{k-1}$ 
along with walks $W_i$, $i=0,1,\ldots{},k-1$, where $W_i$ is a $(v_i,v_{i+2})$-walk
of length $\sigma(i)$ which avoids $v_{i+1}$, such that the total length 
$\sum_{i=0}^{k-1} \sigma(i)$ is odd (subscripts are modulo $k$). 
Kaleidoscopes of order 6 and 7 are depicted in Figure~\ref{kaleidoscopes}.

Let $K_1+K_3$ be the disjoint union of a $K_1: v_0$ and a $K_3: v_1xy$. Then
$K_1+K_3$ contains a kaleidoscope of order 2. Indeed, vertices $v_0, v_1$ along
with the $(v_0,v_0)$-walk $W_0: v_0$ of length 0 (which avoids $v_1$) and 
the $(v_1,v_1)$-walk $W_1: v_1xyv_1$ of length 3 (which avoids $v_0$)
form a kaleidoscope.

Consider an odd cycle $C_{2q+1}: v_0v_1 \dots v_{2q}$ where $q \geq 2$. It contains 
a kaleidoscope consisting of $v_0,v_q,v_{2q}, \dots, v_{q+1}$ along with walks 
$W_i: v_iv_{i-1}$ (which avoids $v_{i+q}$) where $0 \leq i \leq 2q$ 
(subscripts are modulo $2q+1$).
The complement $\overline{C_{2q+1}}$ of $C_{2q+1}$ also contains a kaleidoscope
consisting of $v_0, v_1, \dots, v_{2q}$ along with walks $W_i: v_iv_{i+2}$ (which
avoids $v_{i+1}$) where $0 \leq i \leq 2q$. Thus every odd cycle of length $\geq 5$ 
and its complement contain kaleidoscopes.

\begin{prp} \label{kal2}
If $G$ contains a kaleidoscope of order 2, then it contain $K_1+K_3$ or an odd
cycle of length $\geq 5$ as an induced subgraph.
\end{prp}
\pf Let a kaleidoscope of order 2 in $G$ consist of vertices $v_0,v_1$ along with
walks $W_0$ and $W_1$ of lengths $\sigma(0)$ and $\sigma(1)$ respectively where 
$\sigma(0)+\sigma(1)$ is odd.
Since $\sigma(0)+\sigma(1)$ is odd, precisely one of $\sigma(0)$ and $\sigma(1)$ 
is odd. Without loss of generality assume $\sigma(0)$ is odd. Then $W_0$ is 
a $(v_0,v_0)$-walk of odd length which avoids $v_1$. Since $W_0$ is a closed walk
of odd length, it contains an induced odd cycle $C$. Either $C$ is an induced odd 
cycle of length $\geq 5$ in $G$ or $C$ is a $K_3$ in which case $v_1$ and $C$ induce
a $K_1+K_3$ in $G$. 
\qed

\medskip

Gallai \cite{gallai} introduced the concept of asteroids in the study of
comparability graphs. An {\bf asteroid} of order $k$ in a graph $G$ where $k \geq 3$
is odd consists of vertices $v_0, v_1, \dots, v_{k-1}$ along with walks $W_i$, 
$i = 0, 1, \dots, k-1$, where $W_i$ is a $(v_i,v_{i+2})$-walk that avoids 
$v_{i+1}$. Thus, when the total length of the walks $W_i$ is odd, an asteroid 
is a kaleidoscope. An orientation of a graph $G$ is {\bf transitive}
(respectively, {\bf quasi-transitive}) if for any three vertices $x, y, z$, if
$xy, yz$ are arcs then $xz$ is an arc (respectively, either $xz$ or $zx$ is an arc)
\cite{bh}. Clearly, a transitive orientation is a quasi-transitive orientation but
the converse is not true. Nevertheless, the underlying graphs of transitive
orientable graphs coincide with the underlying graphs of quasi-transitive orientable
graphs \cite{ghouri,hh}. They are called {\bf comparability} graphs.
Gallai \cite{gallai} proved that a graph is a comparability graph if and only if its
complement does not contain an asteroid. An analogous statement holds for strong
comparability graphs which are a subclass of comparability graphs \cite{huangGallai}.
We prove yet another analogous statement for locally complete 2-edge-colourable 
graphs.

\begin{tm} \label{kal}
A graph $G$ is locally complete 2-edge-colourable if and only if its complement
$\overline{G}$ does not contain a kaleidoscope.
\end{tm}
\pf Suppose that $\overline{G}$ contains a kaleidoscope of order $k$. Let
$v_0, v_1, \dots, v_{k-1}$ along with walks $W_i$ connecting $v_i$ and $v_{i+2}$
be a kaleidoscope in $\overline{G}$.
Denote $W_i: v_{i,0}v_{i,1}\dots v_{i,\sigma(i)}$ where $v_{i,0} = v_i$ and
$v_{i,\sigma(i)} = v_{i+2}$, and $\sigma(i) = \sigma(W_i)$ for each $i$.
Since $W_i$ avoids $v_{i+1}$ in $\overline{G}$, $v_{i+1}v_{i,j}$ and 
$v_{i+1}v_{i,j+1}$ form an induced path in $G$ and hence they are adjacent vertices 
in $G^+$ for each $j = 0, 1, \dots, \sigma(i)-1$. This means that
\[Z_{i+1}:\ v_{i+1}v_{i,0}, v_{i+1}v_{i,1}, \dots, v_{i+1}v_{i,\sigma(i)}\]
is a walk in $G^+$ of length $\sigma(i)$ for each $i = 0, 1, \dots, k-1$.
Concatenating these walks $Z_i$ ($0 \leq i \leq k-1$) we obtain a closed walk in
$G^+$ of length $\sum_{i=0}^{k-1} \sigma(i)$ which is odd. Hence $G^+$ is not
bipartite and by Theorem \ref{aux} $G$ is not locally complete 2-edge-colourable.

Conversely, suppose that $G$ is not locally complete 2-edge-colourable. Then
$G^+$ contains an odd cycle according to Theorem \ref{aux}.
Let $e_0e_1 \dots e_{2p}$ be an odd cycle in $G^+$. By definition $e_i$ and
$e_{i+1}$ are edges of an induced path in $G$ and in particular they share an
end vertex for each $i = 0, 1, \dots, 2p$ where subscripts are modulo $2p+1$.
Denote $e_i = u_iv_i$ for each $i = 0, 1, \dots, 2p$. Suppose that
$u_0 = u_1 = \cdots = u_{2p} = u$ (i.e., $u$ is the common end vertex of all edges
$e_i$). Since $e_i$ and $e_{i+1}$ form an induced path in $G$, $v_i$ is
not adjacent to $v_{i+1}$ in $G$ for each $i$. Thus $v_0v_1 \dots v_{2p}$ is an odd
closed walk in $\overline{G}$ and hence it must contain an induced odd cycle $C$.
Either $C$ has length $\geq 5$ or is a $K_3$ in which case $C+u$ is an induced
$K_1+K_3$ in $\overline{G}$. In either case we know from the above that 
$\overline{G}$ contains a kaleidoscope. So we may assume that no vertex is the 
common end vertex of all edges $e_i$.

Call $e_i$ {\bf distinguished} if the common end vertex of $e_{i-1}$ and $e_i$ is
distinct from the common end vertex of $e_i$ and $e_{i+1}$.
Since no vertex is the common end vertex of all edges $e_i$, there is at least one
distinguished edge. Let $e_{i_0}, e_{i_1}, \dots, e_{i_{k-1}}$ where
$0 \leq i_1 < i_2, \dots, < i_k \leq 2p$ be the distinguished edges,
that is, if $i_j < i < i_{j+1}$ then $e_i$ is not distinguished for each
$j = 0, 1, \dots, k-1$.
This means that the edges $e_{i_j}, e_{i_j+1}, \dots, e_{i_{j+1}-1}, e_{i_{j+1}}$
share the same end vertex. Without loss of generality assume that their shared
end vertex is $v_{i_j} = u_{i_j+1} = \cdots = u_{i_{j+1}}$ for each 
$j = 0, 1, \dots, k-1$.
Then $W_j: v_{i_{j-1}}v_{i_j+1}v_{i_j+2}\dots v_{i_{j+1}-1}v_{i_{j+1}}$ is a walk
in $\overline{G}$ connecting $v_{i_{j-1}}$ and $v_{i_{j+1}}$ which avoids
$v_{i_j}$ for each $j=0, 1, \dots, k-1$.
The length $\sigma(j) = i_{j+1}-i_j$ for each $j = 0, 1, \dots, k-2$ and
$\sigma(k-1) = 2p+1+i_0-i_{k-1}$. So the total length
$\sum_{j=0}^{k-1} \sigma(j)= 2p+1$ is odd.
Therefore the vertices $v_{i_0}, v_{i_1}, \dots, v_{i_{k-1}}$ along
with the walks $W_j$ form a kaleidoscope of order $k$ in $\overline{G}$.
\qed

\section{Graphs which have locally complete 2-edge-colourings but are not proper 
         circular-arc graphs}\label{sec3}

The similarity between locally complete 2-edge-coloured graphs and local 
tournaments prompts the following questions: which graphs have locally complete 
2-edge-colourings but cannot be oriented as local tournaments, and which graphs can 
be oriented as local tournaments but do not have locally complete 2-edge-colourings.

The former can be answered by combining a theorem of Skrien \cite{skrien} and a
theorem of Tucker \cite{tucker}. 

\begin{tm} \label{skrien} \cite{skrien}
A connected graph has a local tournament orientation if and only if it is a proper
circular-arc graph.
\qed
\end{tm}

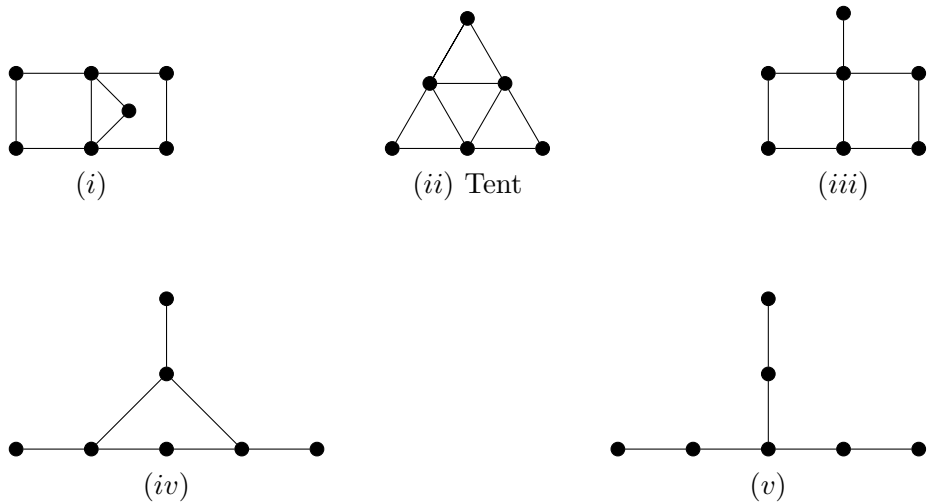
\begin{figure}[!h]
\begin{center}
\begin{tikzpicture}[>=latex]
%\begin{pgfonlayer}{nodelayer}
                \node [style=blackvertex] (1) at (6, 1.73) {};
                \node [style=blackvertex] (2) at (5.5, 0.865) {};
                \node [style=blackvertex] (3) at (6.5, 0.865) {};
                \node [style=blackvertex] (4) at (5, 0) {}; 
                \node [style=blackvertex] (5) at (6, 0) {}; 
                \node [style=blackvertex] (6) at (7, 0) {}; 

\draw [-] (2) -- (3) -- (5) -- (2) -- (1) -- (3) -- (6) -- (5) -- (4)-- (2) -- (1);
\node [style=textbox]  at (1, -.5) {$(i)$};

                        \node[unlabelled] (7) at (11, 1.8)       {};
                        \node[unlabelled] (8) at (10, 1) {};
                        \node[unlabelled] (9) at (11, 1) {};
                        \node[unlabelled] (10) at (12, 1) {};
                        \node[unlabelled] (11) at (10, 0) {};
                        \node[unlabelled] (12) at (11, 0) {};
                        \node[unlabelled] (13) at (12, 0) {};

\draw [-] (7) -- (9) -- (12) -- (13) -- (10) -- (9) -- (8) -- (11) -- (12); 

\node [style=textbox]  at (6, -.5) {$(ii)$\ Tent};

                        \node[unlabelled] (14) at (0, 1) {};
                        \node[unlabelled] (15) at (1, 1) {};
                        \node[unlabelled] (16) at (2, 1) {};
                        \node[unlabelled] (17) at (0, 0) {};
                        \node[unlabelled] (18) at (1, 0) {};
                        \node[unlabelled] (19) at (2, 0) {};
                        \node[unlabelled] (20) at (1.5, 0.5)     {};

\draw [-] (15) -- (20) -- (18) -- (15) -- (16) -- (19) -- (18) -- (17) -- (14) -- 
             (15);

\node [style=textbox]  at (11, -.5) {$(iii)$};

                        \node[unlabelled] (21) at (0, -4) {};
                        \node[unlabelled] (22) at (1, -4) {};
                        \node[unlabelled] (23) at (2, -4) {};
                        \node[unlabelled] (24) at (3, -4) {};
                        \node[unlabelled] (25) at (4, -4) {};
                        \node[unlabelled] (26) at (2, -3) {};
                        \node[unlabelled] (27) at (2, -2) {};

\draw [-] (21) -- (22) -- (23) -- (24) -- (25);
\draw [-] (22) -- (26) -- (27); 
\draw [-] (26) -- (24);
\node [style=textbox]  at (2, -4.5) {$(iv)$};

                       \node[unlabelled] (28) at (8, -4) {};
                        \node[unlabelled] (29) at (9, -4) {};
                        \node[unlabelled] (30) at (10, -4) {};
                        \node[unlabelled] (31) at (11, -4) {};
                        \node[unlabelled] (32) at (12, -4) {};
                        \node[unlabelled] (33) at (10, -3) {};
                        \node[unlabelled] (34) at (10,-2) {};
\draw [-] (28) -- (29) -- (30) -- (31) -- (32);
\draw [-] (34) -- (33) -- (30);

\node [style=textbox]  at (10, -4.5) {$(v)$};
                
%\end{pgfonlayer}
\end{tikzpicture}
\end{center}
\vspace{-2mm}
\caption{Some forbidden subgraphs for the complements of proper circular-arc graphs.}
\label{tuckerlist}
\end{figure}

\begin{tm} \label{tucker} \cite{tucker}
A connected graph $G$ is a proper circular-arc graph if and only if $G$ does not contain
$C_k + K_1$ ($k \geq 4$) or tent $+ K_1$ and $\overline{G}$ does not contain
$C_{2k}$ ($k \geq 3$), $C_{2k+1}+K_1$ ($k \geq 1$), or any of the graphs in
Figure~\ref{tuckerlist} as an induced subgraph.
\qed
\end{tm}

Suppose that $G$ is a locally complete 2-edge-colourable graph which cannot be 
oriented as a local tournament. Since $G$ is locally complete 2-edge-colourable,
every component of $G$ is also locally complete 2-edge-colourable. Since $G$ cannot
be oriented as a local tournament, at least one component of $G$ cannot be oriented
as a local tournament. Let $H$ be such a component. Thus $H$ is a connected graph
which is locally complete 2-edge-colourable but cannot be oriented as a local 
tournament. By Theorems \ref{skrien} and \ref{tucker}, $H$ contains 
$C_k + K_1$ ($k \geq 4$) or tent $+ K_1$, or $\overline{H}$ contains
$C_{2k}$ ($k \geq 3$), $C_{2k+1}+K_1$ ($k \geq 1$), or any of the graphs in
Figure~\ref{tuckerlist} as an induced subgraph. We know that neither $C_{2k+1}$ nor 
$\overline{C_{2k+1}}$ has a complete locally complete 2-edge-colouring for each 
$k \geq 2$. It is also easy to check that neither Tent (Figure~\ref{tuckerlist}$(ii)$) nor
the complement of Figure~\ref{tuckerlist}$(i)$ has a locally complete 
2-edge-colouring and that the complements of each of the graphs in Figure \ref{tuckerlist} (ii)-(v) have such a colouring. Hence we have the following:

\begin{tm} \label{lc2ec-minus-pca}
Let $G$ be a locally complete 2-edge-colourable graph. Then $G$ is not a proper 
circular-arc graph (i.e., cannot be oriented as a local tournament) if and only if
it contains $C_{2k}+K_1$ for some $k \geq 2$, $\overline{C_{2k}}$ for some 
$k \geq 3$, or the complement of a graph in Figure~\ref{tuckerlist}$(ii)-(v)$ as
an induced subgraph.
\qed
\end{tm}

\section{Proper circular-arc graphs which are locally 2-edge-colourable}\label{sec4}

In this section we study the graphs which can be oriented as local tournaments and 
also have locally complete 2-edge-colourings. For this purpose we can restrict 
our study to connected graphs and hence by Theorem \ref{skrien} to connected 
proper circular arc graphs. 

Two vertices in a graph are {\bf true-twins} if they have the same closed 
neighbourhood. A graph is {\bf reduced} if it does not have true twins. 
Given a graph $G$, a reduced graph $G'$ can be obtained from $G$ by repeatedly 
deleting one of true twins until no two vertices are true twins. Clearly, such
a reduced graph $G'$ is unique and we call $G'$ the {\bf reduced subgraph} of $G$.
One verifies easily that $G$ can be oriented as a local tournament if and only 
if its reduced subgraph can be oriented as a local tournament, and that $G$ has 
a locally complete 2-edge-colouring if and only if its reduced subgraph has 
a locally complete 2-edge-colouring. Thus we can further restrict our study to 
reduced graphs. 

A {\bf universal} vertex in a graph $G$ is a vertex adjacent to every other vertex.
A reduced graph has at most one universal vertex.

Let $G$ be a connected reduced proper circular-arc graph. Two edges in $G$ are in 
the same {\bf implication class} if they form an induced path (of length 2) in $G$.
It follows from definition of a locally complete 2-edge-colouring that edges in the same implication class must have 
different colours in a locally complete 2-edge-colouring of $G$. 

It is proved in \cite{huang} that if $\overline{G}$ is not bipartite or 
$\overline{G}$ is bipartite and $G$ is a proper interval graph without a universal 
vertex then all edges of $G$ are in the same implication class. Thus we have 
the following:

\begin{tm} \label{all-in-one}
Let $G$ be a connected reduced proper circular-arc graph which is locally complete 
2-edge-colourable. If $\overline{G}$ is not bipartite or $\overline{G}$ is bipartite 
and $G$ is a proper interval graph without universal vertex, then $G$ has a locally
complete 2-edge-colouring, which is unique up to switching the colours. 
\qed
\end{tm}
 
A {\bf straight ordering} of a graph $G$ is a linear ordering of the vertices of $G$,
$v_1, v_2, \dots, v_n$, such that for each $i$,
\[N[v_i] = \{v_{i-\ell}, v_{i-\ell+1}, \dots, v_{i-1}, v_i, 
             v_{i+1}, v_{i+2}, \dots, v_{i+\gamma}\}. \]
Here $\ell$ and $\gamma$ are non-negative and depend on $i$.
Similarly, a {\bf round ordering} of $G$ is a circular ordering of the 
vertices of $G$, $v_1, v_2, \dots, v_n$, such that for each $i$,
\[N[v_i] = \{v_{i-\ell}, v_{i-\ell+1}, \dots, v_{i-1}, v_i, 
             v_{i+1}, v_{i+2}, \dots, v_{i+\gamma}\} \]
where $\ell$ and $\gamma$ are non-negative and depend on $i$, and the subscript
additions and subtractions are modulo $n$. 
If needed, we shall write $\ell(i)$ and $\gamma(i)$ instead of just $\ell$ and 
$\gamma$ in either case of a straight ordering or a round ordering. For convenience, 
we shall also let $\ell(v_i) = v_{i-\ell(i)}$ and $\gamma(v_i) = v_{i+\gamma(i)}$.     
\begin{prp} \label{sr-ordering} \cite{huang}
The following statements hold for any graph $G$:
\begin{itemize}
\item[(1)] $G$ is a proper interval graph if and only if it has a straight ordering.
\item[(2)] $G$ is a proper circular-arc graph if and only if it has a round ordering.
\qed
\end{itemize}
\end{prp}

  \begin{remark}
    \label{pig/pcastructure}
  
When $v_1, v_2, \dots, v_n$ is either a straight or a round ordering of $G$, we 
shall use $G[v_i,v_j]$ for any $i \leq j$ to denote the subgraph of $G$ induced by 
$v_i, v_{i+1}, \dots, v_j$ and $G[v_j,v_i]$ to denote the subgraph of $G$ induced by 
$v_j, v_{j+1}, \dots, v_n, v_1, v_2, \dots, v_i$. To simplify notation further,
we also use $G[i,j]$ instead of $G[v_i,v_j]$ and $G[j,i]$ instead of $G[v_j,v_i]$. 
By results in \cite{huang} we can assume that
\begin{equation}\label{complete}
  G[\ell(i),i] \mbox{ and } G[i,\gamma(i)] \mbox{ are 
    complete subgraphs of } G.
\end{equation}

\end{remark}

A {\bf clique covering} of a graph $G$ is a set of cliques in $G$ that covers $G$
in the way that each vertex of $G$ is in at least one of the cliques.
The {\bf clique covering number} of $G$, denoted by $cc(G)$ is the least
integer $k$ such that there is a clique covering of $G$ with $k$ cliques.
When $cc(G) \leq 2$, $\overline{G}$ is bipartite and,  as observed in Section~1,  
$G$ has a locally complete 2-edge-colouring. We state this observation for easy 
reference in the rest of the paper.

\begin{prp} \label{cc2}
If $cc(G) \leq 2$, then $G$ is locally complete 2-edge-colourable.
\qed
\end{prp}

\subsection{The case when $G$ is a proper interval graph}

Our goal in this subsection is to give a full description of (connected reduced) 
proper interval graphs which have locally complete 2-edge-colourings.

\begin{lemma} \label{lem:cc=2}
Let $G$ be a connected reduced graph with $cc(G) = 2$ and without a cutvertex. 
Suppose that $G$ is a proper interval graph and $v_1,v_2,\ldots{},v_n$ is a straight
ordering of $G$. Then the following hold:
\begin{itemize}
\item[(1)] If $G$ has no universal vertex, then for every locally complete
       2-edge-colouring of $G$, all edges incident with at least one of $v_1$ and $v_n$ have 
       the same colour. 
\item[(2)] If $G$ has a universal vertex, then for every locally complete 
       2-edge-colouring of $G$, exactly one of $v_1$ and $v_n$ is incident with 
       edges of both colours.
\end{itemize}
\end{lemma}
\pf Suppose first that $G$ has no universal vertex. As $cc(G)=2$, there is a vertex 
$v_i$ such that $v_1\sim v_i \nad v_n$ and $v_1 \nad v_{i+1}\sim v_n$.
We colour all edges of $G[1,i]$ and all edges of $G[i+1,n]$ by colour 1 and all 
remaining edges by colour 2. This is a locally complete 2-edge-colouring in which
all edges incident with $v_1$ and $v_n$ have the same colour. By
Theorem \ref{all-in-one} this colouring is unique up to switching the colours. 
Hence (1) holds.

Suppose that $v_c$ is the universal vertex. As $cc(G) = 2$, $1 < c < n$ an thus
$v_1 \sim v_c \sim v_n$. Since $G$ is reduced, $\gamma(1) = c = \ell(n)$,
$\gamma(v_2) = v_{c+1}$ and $\ell(v_{n-1}) = v_{c-1}$. In particular, $v_1v_cv_n$, 
$v_1v_iv_{c+1}v_n$ and $v_1v_{c-1}v_jv_n$ are induced paths in $G$ for each 
$1 < i < c$ and $c < j < n$. Consider an arbitrary locally complete 2-edge-colouring 
$\phi: E(G) \rightarrow \{1,2\}$ of $G$. Assume without loss of generality that 
$\phi(v_1v_c)=1$. Since $v_1v_cv_n$ is an induced path, $\phi(v_cv_n)=2$. 
If $\phi(v_{c-1}v_{c+1}) = 1$, then $\phi(v_1v_{c-1}) = 2$ as $v_1v_{c-1}v_{c+1}$ is 
an induced path. Since each $v_1v_{c-1}v_jv_n$ with $c < j < n$ is an induced path,
$\phi(v_{c-1}v_j) = 1$ and $\phi(v_jv_n) = 2$ for each $c < j < n$. Hence 
edges incident with $v_1$ have both colours and the edges incident with $v_n$ 
have only colour 2. Similarly, if $\phi(v_{c-1}v_{c+1}) = 2$, then 
$\phi(v_{c+1}v_n) = 1$ as $v_{c-1}v_{c+1}v_n$ is an induced path. 
Since each $v_1v_iv_{c+1}v_n$ with $1 < i < c$ is an induced path, 
$\phi(v_1v_i) = 1$ for each $1 < i < c$. Hence the edges incident with $v_1$ have 
only colour 1 and the edges incident with $v_n$ have both colours.  
\qed

\begin{figure}[ht]
\begin{center}
\begin{tikzpicture}[>=latex]
%\begin{pgfonlayer}{nodelayer}
                \node [label={below: 1}] [style=blackvertex] (1) at (1,0) {};
                \node [label={below: 2}] [style=blackvertex] (2) at (2,0) {};
                \node [label={below: 3}] [style=blackvertex] (3) at (3,0) {};
                \node [label={below: 4}] [style=blackvertex] (4) at (4,0) {};
                \node [label={below: 5}] [style=blackvertex] (5) at (5,0) {};
                \node [label={below: 6}] [style=blackvertex] (6) at (6,0) {};
                \node [label={below: 7}] [style=blackvertex] (7) at (7,0) {};

\draw[-] (1) -- (2) -- (3) -- (4) -- (5) -- (6) -- (7);
\draw [ ] (2) to [out=30,in=150] (4);
\draw [ ] (3) to [out=30,in=150] (5);
\draw [ ] (4) to [out=30,in=150] (6);

\node [style=textbox]  at (9, 0) {$F_1$};

                \node [label={below: 1}] [style=blackvertex] (8) at (1,-2) {};
                \node [label={below: 2}] [style=blackvertex] (9) at (2,-2) {};
                \node [label={below: 3}] [style=blackvertex] (10) at (3,-2) {};
                \node [label={below: 4}] [style=blackvertex] (11) at (4,-2) {};
                \node [label={below: 5}] [style=blackvertex] (12) at (5,-2) {};
                \node [label={below: 6}] [style=blackvertex] (13) at (6,-2) {};
                \node [label={below: 7}] [style=blackvertex] (14) at (7,-2) {};

\draw[-] (8) -- (9) -- (10) -- (11) -- (12) -- (13) -- (14);
\draw [ ] (8) to [out=30,in=150] (10);
\draw [ ] (9) to [out=30,in=150] (11);
\draw [ ] (10) to [out=30,in=150] (12);
\draw [ ] (11) to [out=30,in=150] (13);
\draw [ ] (12) to [out=30,in=150] (14);

\node [style=textbox]  at (9, -2) {$F_2$};

                \node [label={below: 1}] [style=blackvertex] (15) at (1,-4) {};
                \node [label={below: 2}] [style=blackvertex] (16) at (2,-4) {};
                \node [label={below: 3}] [style=blackvertex] (17) at (3,-4) {};
                \node [label={below: 4}] [style=blackvertex] (18) at (4,-4) {};
                \node [label={below: 5}] [style=blackvertex] (19) at (5,-4) {};
                \node [label={below: 6}] [style=blackvertex] (20) at (6,-4) {};
                \node [label={below: 7}] [style=blackvertex] (21) at (7,-4) {};

\draw[-] (15) -- (16) -- (17) -- (18) -- (19) -- (20) -- (21);
\draw [ ] (15) to [out=30,in=150] (17);
\draw [ ] (16) to [out=30,in=150] (18);
\draw [ ] (17) to [out=30,in=150] (19);
\draw [ ] (18) to [out=30,in=150] (20);
\draw [ ] (19) to [out=30,in=150] (21);
\draw [ ] (16) to [out=35,in=145] (19);
\draw [ ] (17) to [out=35,in=145] (20);

\node [style=textbox]  at (9, -4) {$F_3$};

%\end{pgfonlayer}
\end{tikzpicture}
\end{center}
\vspace{-2mm}
\caption{Proper interval graphs $F_1$, $F_2$ and $F_3$.}
\label{fs}
\end{figure}
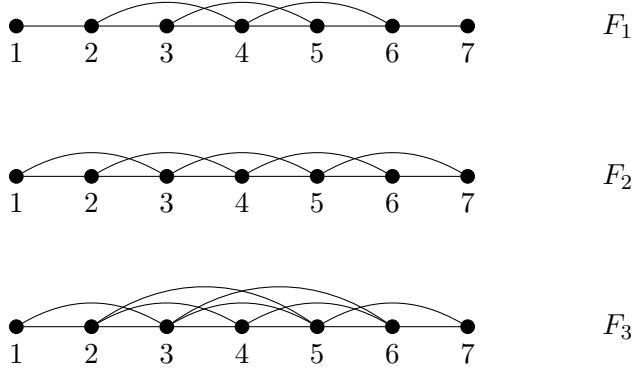

\begin{prp} \label{prop:Fhasnocol}
None of the graphs $F_1,F_2,F_3$ in Figure~\ref{fs} is locally complete 2-edge-colourable.
\end{prp}
\pf There are odd cycles in the auxiliary graph $F_i^+$ of $F_i$ for each 
$i = 1, 2, 3$. Indeed, $12\sim 24\sim 46\sim 67\sim 56\sim 35\sim 23\sim 12$,
$13\sim 35\sim 57\sim 45\sim 24\sim 46\sim 34\sim 13$ and 
$13\sim 36\sim 67\sim 46\sim 24\sim 12\sim 25\sim 57\sim 35\sim 13$ are 
odd cycles in $F_1^+, F_2^+$ and $F_3^+$ respectively. Hence by Theorem \ref{aux}
none of $F_1,F_2,F_3$ has a locally complete 2-edge-colouring.
  \qed

\medskip

It follows from Proposition \ref{prop:Fhasnocol} that if a graph contains any of 
$F_1, F_2, F_3$ as in induced subgraph, then it does not have a locally complete 2-edge-colouring. The 
straight orderings given in Figure \ref{fs} show that $F_1, F_2, F_3$ are proper
interval graphs. It can be checked that deleting any vertex of these graphs results
in a locally complete 2-edge-colourable graph. Hence $F_1, F_2, F_3$ are (minimal)
forbidden subgraphs for the class of proper interval graphs which have locally 
complete 2-edge-colourings. 

Let $v_1, v_2, \dots, v_n$ be a straight ordering of a proper interval graph $G$.
Then for each edge $v_iv_j$ with $i < j$, $G[i,j]$ is a clique. This implies that each 
maximal clique in $G$ is of form $G[i,j]$ for some $i \leq j$.

Let $k=cc(G)$. The {\bf canonical clique covering} of $G$ with respect to 
the straight ordering $v_1, v_2, \dots, v_n$ consists of maximal cliques 
$G[s_i,t_i]$ ($1 \leq i \leq k$) such that $t_i$ is the largest integer with 
$cc(G[1,t_i]) = i$ for each $i = 1, 2, \dots, k$. We summarize some basic properties
of the canonical clique covering  in the following lemma.

\begin{lemma} \label{canonical}
Let $G$ be a connected proper interval graph with $cc(G) = k$ and let 
$v_1, v_2, \dots, v_n$ be a straight ordering of $G$. Suppose that $G[s_i,t_i]$ 
($1 \leq i \leq k$) are the cliques of the canonical clique covering of $G$. Then
the following hold:
\begin{itemize}
\item[(1)] $s_1 = 1$ and $t_k = n$;
\item[(2)] $\gamma(v_{s_i}) = v_{t_i}$ and $\ell(v_{t_i}) = v_{s_i}$ for each 
           $i = 1, 2, \dots, k$;
\item[(3)] $s_1 < s_2 < \cdots < s_k$ and $t_1 < t_2 < \cdots < t_k$;
\item[(4)] if $G$ is reduced then $s_i = t_{i-1}$ or $s_i=t_{i-1}+1$ for each $i=2,3,\dots{},k$.
\end{itemize}
\end{lemma} 
\pf Since each $G[s_i,t_i]$ is a maximal clique, $\gamma(v_{s_i}) = v_{t_i}$ and 
$\ell(v_{t_i}) = v_{s_i}$ for each $i$. This proves $(2)$. Since $t_i$ is largest 
integer with $cc(G[1,t_i]) = i$ for each $i$, it follows that a
$t_1 < t_2 < \cdots < t_k$. Combining
this with $(2)$ we have $s_1 < s_2 < \cdots < s_k$. Thus $(3)$ holds and $(1)$ 
follows from it. We must have $s_i \leq t_{i-1}+1$ as otherwise $v_{t_{i-1}+1}$ does
not belong to any of the cliques $G[s_i,t_i]$, a contradiction to the assumption.
If $s_i < t_{i-1}$, then $N[v_{s_i}] = N[v_{t_{i-1}}] = V(G[s_{i-1},t_i])$, which 
means that $G$ is not reduced. This proves $(4)$.
\qed

\begin{lemma} \label{2-canonical}
Let $G$ be a 2-connected reduced proper interval graph with $cc(G) = k \geq 2$.
Let  $v_1,v_2,\ldots{},v_n$ be a straigh ordering of $V(G)$ and 
let $G[s_i,t_i]$ ($1 \leq i \leq k$) be the canonical clique covering of $G$ .
 Then $G$ contains 
one of the graphs $F_1, F_2, F_3$ in Figure \ref{fs} as an induced subgraph unless
the following hold:
\begin{itemize}
\item[(1)] For each $i$ with $1 < i < k$, if $\gamma(v_{s_i-1}) = v_a$ and 
  $\ell(v_{t_i+1}) = v_b$ then $a < b$.
  \item[(2)] If $t_i = s_{i+1}$, then $i = 1$ or $i=k-1$.
\end{itemize}
\end{lemma}
\pf To prove that $(1)$ holds suppose for a contradiction that for some $i$ with $1<i<k$ we have $a \geq b$. Let $1<i<k$ be chosen as small as possible such that $\gamma(v_{s_i-1})\geq \ell(v_{t_i+1})$. By Lemma \ref{canonical} $(4)$, 
$s_i = t_{i-1}$ or $t_{i-1}+1$, and $t_i = s_{i+1}$ or $s_{i+1}-1$.\\ 

\noindent{}{\bf Case~1.} $s_i = t_{i-1}$ and $t_i = s_{i+1}$.

Since $a \geq b$ and $v_b \sim v_{t_i+1}$, $v_a \sim v_{t_i+1} = v_{s_{i+1}+1}$. 
We see that $G$ contains a copy of $F_2$ induced by 
$v_{s_{i-1}}, v_{s_i-1}, v_{s_i}, v_a, v_{t_i} v_{t_i+1}, t_{i+1}$.\\

\noindent{}{\bf Case 2.} $s_i = t_{i-1}$ and $t_i = s_{i+1}-1$.

Since $v_b = \ell(v_{t_i+1})$, $G[b,t_i+1]$ is a clique. We must have $b > s_i+1$ as otherwise $b = s_i+1 = t_{i-1}+1$ and
$cc(G[1,t_i+1] = cc(G[1,t_i])$, a contradiction
to the definition of the canonical clique covering. Hence $G$ contains a copy of 
$F_3$ induced by $v_{s_{i-1}}, v_{s_i-1}, v_{s_i}, v_{s_i+1}, v_a, v_{t_i}, 
v_{t_i+1}$.\\

\noindent{}{\bf Case 3.} $s_i = t_{i-1}+1$ and $t_i = s_{i+1}$.

Since $G$ has no cutvertex, $t_{i-1}$ is not a cutvertex, which means that 
$\ell(v_{s_i})$ is a vertex in $G[s_{i-1},t_{i-1}]$ distinct from $v_{t_{i-1}}$. 
If $\ell(v_{s_i}) \not\sim v_a$, then $v_{s_{i-1}}, \ell(v_{s_i}), v_{t_{i-1}},
v_{s_i}, v_a, v_{t_i}, v_{t_i+1}$ induce a copy of $F_2$.
Otherwise, we have $\ell(v_{s_i}) \sim v_a$ which implies 
$N[v_{t_{i-1}}] \subseteq N[\ell(v_{s_i})]$. Since $G$ is reduced, there must be a
vertex adjacent to $\ell(v_{s_i})$ but not to $v_{t_{i-1}}$. Let $v_j$ be such
a vertex. Then $j < s_{i-1}$. So $i > 2$. Now we see that $v_{s_{i-1}-1}\sim v_{\ell(s_i)}\sim v_{t_{i-1}+1}=v_{s_i}$, so $\gamma(v_{s_{i-1}-1})\geq \ell(v_{t_{i-1}+1})$, contradicting the minimality of $i$.\\

\noindent{}{\bf Case 4.} $s_i = t_{i-1}+1$ and $t_i = s_{i+1}-1$.

We see immediately that $v_{s_{i-1}}, v_{t_{i-1}}, v_{s_i}, v_a, v_{t_i}, 
v_{s_{i+1}}, v_{t_{i+1}}$ induce a copy of $F_1$.\\

To prove $(2)$ suppose that $t_i = s_{i+1}$ for some $i$ with $1 < i < k-1$. Then $k \geq 4$.
Since $G$ has no cut-vertex, each clique in the canonical covering of $G$ must have 
at least 3 vertices. In particular, $t_i-1 > s_i$. Moreover, if 
$\gamma(v_{t_i-1}) = v_j$ then $t_i < j < t_{i+1}$ by Lemma \ref{canonical} (2) and the fact that $G$ has no cutvertex. We see now that
$\ell(v_{s_i}), v_{s_i}, v_{t_i-1}, v_{t_i}, \gamma(v_{t_i-1}), v_{t_{i+1}}, \gamma(v_{t_{i+1}})$
induce a copy of $F_1$ in $G$. Here we used $(1)$ to see that $\ell(v_{s_i})\nad v_{t_i-1}$ and $v_{\gamma(v_{t_i-1})}\nad v_{\gamma(v_{t_{i+1}})}$.
This proves $(2)$. 
\qed

\medskip

By Lemma \ref{2-canonical}(2), if $t_i = s_{i+1}$ and 
$t_{i+1} = s_{i+2}$ then $k \leq 3$. Combining it with Lemma \ref{canonical}(4),
we have that when $k \geq 4$, if $s_2=t_1$ and $t_{k-1}=s_k$ 
%$s_1 = t_2$ and $s_{k-1} = t_k$
then
$s_i = t_{i-1}+1$ for all $3 \leq i \leq k-1$.

\begin{dfn}\label{def:types}
Let $G$ be a 2-connected reduced proper interval graph with $cc(G) = k \geq 2$. % and 
%without cut-vertex. 
Let $G[s_i,t_i]$ ($1 \leq i \leq k$) be the canonical clique covering of $G$ 
with the straight ordering $v_1,v_2,\ldots{},v_n$. We say that $G$ is of
\begin{itemize}
\item {\bf type 1} if $t_1 = s_2$, $t_{k-1} = s_k$ and $t_i = s_{i+1}-1$ 
                   for all $2 \leq i \leq k-2$.
\item {\bf type 2} if $k \geq 3$, $t_1 = s_2$ and $t_i = s_{i+1}-1$ for all 
                   $2 \leq i \leq k-1$.
\item {\bf type 3} if $k \geq 3$, $t_{k-1} = s_k$ and $t_i = s_{i+1}-1$ for all 
                   $1 \leq i \leq k-2$.
\item {\bf type 4} if $t_i  = s_{i+1}-1$ for all $1 \leq i \leq k-1$.
\end{itemize}
In all types, for each $i$ with $2 \leq i \leq k-1$, if $\gamma(v_{s_i-1}) = v_a$ and
          $\ell(v_{t_i+1}) = v_b$ then $a < b$. 
\end{dfn}

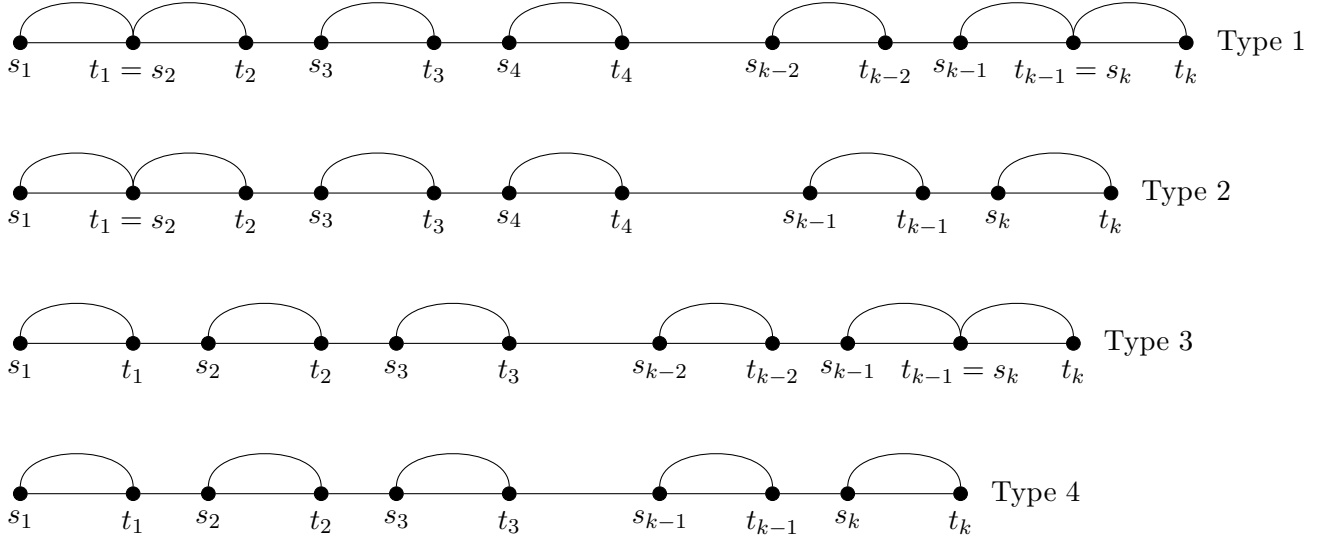
\begin{figure}[ht]
\begin{center}
\begin{tikzpicture}[>=latex]
%\begin{pgfonlayer}{nodelayer}
       \node [label={below: $s_1$}] [style=blackvertex] (1) at (0,0) {};
       \node [label={below: $t_1=s_2$}] [style=blackvertex] (2) at (1.5,0) {};
       \node [label={below: $t_2$}] [style=blackvertex] (3) at (3,0) {};
       \node [label={below: $s_3$}] [style=blackvertex] (4) at (4,0) {};
       \node [label={below: $t_3$}] [style=blackvertex] (5) at (5.5,0) {};
       \node [label={below: $s_4$}] [style=blackvertex] (6) at (6.5,0) {};
       \node [label={below: $t_4$}] [style=blackvertex] (7) at (8,0) {};
       \node [label={below: $s_{k-2}$}] [style=blackvertex] (8) at (10,0) {};
       \node [label={below: $t_{k-2}$}] [style=blackvertex] (9) at (11.5,0) {};
       \node [label={below: $s_{k-1}$}] [style=blackvertex] (10) at (12.5,0) {};
       \node [label={below: $t_{k-1}=s_k$}] [style=blackvertex] (11) at (14,0) {};
       \node [label={below: $t_k$}] [style=blackvertex] (12) at (15.5,0) {};

\draw[-] (1) -- (2) -- (3) -- (4) -- (5) -- (6) -- (7) -- (8) -- (9) -- (10) --
          (11) -- (12);

\draw [ ] (1) to [out=90,in=90] (2);
\draw [ ] (2) to [out=90,in=90] (3);
\draw [ ] (4) to [out=90,in=90] (5);
\draw [ ] (6) to [out=90,in=90] (7);
\draw [ ] (8) to [out=90,in=90] (9);
\draw [ ] (10) to [out=90,in=90] (11);
\draw [ ] (11) to [out=90,in=90] (12);

\node [style=textbox]  at (16.5,0) {Type 1};

       \node [label={below: $s_1$}] [style=blackvertex] (13) at (0,-2) {};
       \node [label={below: $t_1=s_2$}] [style=blackvertex] (14) at (1.5,-2) {};
       \node [label={below: $t_2$}] [style=blackvertex] (15) at (3,-2) {};
       \node [label={below: $s_3$}] [style=blackvertex] (16) at (4,-2) {};
       \node [label={below: $t_3$}] [style=blackvertex] (17) at (5.5,-2) {};
       \node [label={below: $s_4$}] [style=blackvertex] (18) at (6.5,-2) {};
       \node [label={below: $t_4$}] [style=blackvertex] (19) at (8,-2) {};
       \node [label={below: $s_{k-1}$}] [style=blackvertex] (20) at (10.5,-2) {};
       \node [label={below: $t_{k-1}$}] [style=blackvertex] (21) at (12,-2) {};
       \node [label={below: $s_k$}] [style=blackvertex] (22) at (13,-2) {};
       \node [label={below: $t_k$}] [style=blackvertex] (23) at (14.5,-2) {};

\draw[-] (13) -- (14) -- (15) -- (16) -- (17) -- (18) -- (19) -- (20) -- (21) -- 
         (22) -- (23);

\draw [ ] (13) to [out=90,in=90] (14);
\draw [ ] (14) to [out=90,in=90] (15);
\draw [ ] (16) to [out=90,in=90] (17);
\draw [ ] (18) to [out=90,in=90] (19);
\draw [ ] (20) to [out=90,in=90] (21);
\draw [ ] (22) to [out=90,in=90] (23);

\node [style=textbox]  at (15.5,-2) {Type 2};

       \node [label={below: $s_1$}] [style=blackvertex] (24) at (0,-4) {};
       \node [label={below: $t_1$}] [style=blackvertex] (25) at (1.5,-4) {};
       \node [label={below: $s_2$}] [style=blackvertex] (26) at (2.5,-4) {};
       \node [label={below: $t_2$}] [style=blackvertex] (27) at (4,-4) {};
       \node [label={below: $s_3$}] [style=blackvertex] (28) at (5,-4) {};
       \node [label={below: $t_3$}] [style=blackvertex] (29) at (6.5,-4) {};
       \node [label={below: $s_{k-2}$}] [style=blackvertex] (30) at (8.5,-4) {};
       \node [label={below: $t_{k-2}$}] [style=blackvertex] (31) at (10,-4) {};
       \node [label={below: $s_{k-1}$}] [style=blackvertex] (32) at (11,-4) {};
       \node [label={below: $t_{k-1}=s_k$}] [style=blackvertex] (33) at (12.5,-4) {};
       \node [label={below: $t_k$}] [style=blackvertex] (34) at (14,-4) {};

\draw[-] (24) -- (25) -- (26) -- (27) -- (28) -- (29) -- (30) -- (31) -- (32) --
         (33) -- (34);

\draw [ ] (24) to [out=90,in=90] (25);
\draw [ ] (26) to [out=90,in=90] (27);
\draw [ ] (28) to [out=90,in=90] (29);
\draw [ ] (30) to [out=90,in=90] (31);
\draw [ ] (32) to [out=90,in=90] (33);
\draw [ ] (33) to [out=90,in=90] (34);

\node [style=textbox]  at (15,-4) {Type 3};

       \node [label={below: $s_1$}] [style=blackvertex] (35) at (0,-6) {};
       \node [label={below: $t_1$}] [style=blackvertex] (36) at (1.5,-6) {};
       \node [label={below: $s_2$}] [style=blackvertex] (37) at (2.5,-6) {};
       \node [label={below: $t_2$}] [style=blackvertex] (38) at (4,-6) {};
       \node [label={below: $s_3$}] [style=blackvertex] (39) at (5,-6) {};
       \node [label={below: $t_3$}] [style=blackvertex] (40) at (6.5,-6) {};
       \node [label={below: $s_{k-1}$}] [style=blackvertex] (41) at (8.5,-6) {};
       \node [label={below: $t_{k-1}$}] [style=blackvertex] (42) at (10,-6) {};
       \node [label={below: $s_k$}] [style=blackvertex] (43) at (11,-6) {};
       \node [label={below: $t_k$}] [style=blackvertex] (44) at (12.5,-6) {};

\draw[-] (35) -- (36) -- (37) -- (38) -- (39) -- (40) -- (41) -- (42) -- (43) --
         (44);
    
\draw [ ] (35) to [out=90,in=90] (36);
\draw [ ] (37) to [out=90,in=90] (38);
\draw [ ] (39) to [out=90,in=90] (40);
\draw [ ] (41) to [out=90,in=90] (42);
\draw [ ] (43) to [out=90,in=90] (44);

\node [style=textbox]  at (13.5,-6) {Type 4};

%\end{pgfonlayer}
\end{tikzpicture}
\end{center}
\vspace{-2mm}
\caption{Proper interval graphs of types 1 - 4.}
\label{pic:types}
\end{figure}

Figure \ref{pic:types} depicts these 4 types of proper interval graphs. 

\begin{cor} \label{types}
  Let $G$ be a 2-connected reduced proper interval graphs with $cc(G) \geq 2$.
  %and without cutvertex. 
Suppose that $G$ does not contain any of $F_1, F_2$ or $F_3$ as an induced subgraph.
Then $G$ is of one of the 4 types in Definition \ref{def:types}.
\end{cor}
\pf This follows from Lemma \ref{2-canonical}.
\qed

\begin{lemma}\label{lem:PIGcc>=3}
Let $G$ be a reduced proper interval graph which is of one of the 4 types in
Definition \ref{def:types}. Then $G$ has a locally complete 2-edge-colouring 
$\phi$. Moreover, the colouring $\phi$ is unique up to switching the colours, 
except when $cc(G) = 2$ and $G$ is of type 2 (or 3). 
\end{lemma}
\pf We define an edge-colouring $\phi(v_pv_q)$ for all edges $v_pv_q$ with $p < q$ as follows: 

When $G$ is of type 1, 
$$\phi{}(v_pv_q)=\left\{\begin{array}{ll}
  1 &\mbox{ if } q < t_1\ \mbox{or}\ p > s_k\ \mbox{or}\
    s_i \leq p < q \leq t_i\ \mbox{for all}\ 2 \leq i \leq k-1\\
  2 & \mbox{ otherwise }\end{array}\right.$$

When $G$ is of type 2, 
$$\phi{}(v_pv_q)=\left\{\begin{array}{ll}
  1 & \mbox{ if }\mbox q < t_1 \mbox{ or }
        s_i \leq p < q \leq t_i\ \mbox{for all}\ 2 \leq i \leq k\\
  2 & \mbox{ otherwise }\end{array}\right.$$
           
When $G$ is of type 3,
$$\phi{}(v_pv_q)=\left\{\begin{array}{ll}
  1 & \mbox{ if }\mbox p > s_k \mbox{ or }
        s_i \leq p < q \leq t_i\ \mbox{for all}\ 1 \leq i \leq k-1\\
  2 & \mbox{ otherwise }\end{array}\right.$$

When $G$ is of type 4, 
$$\phi{}(v_pv_q)=\left\{\begin{array}{ll}
  1 & \mbox{ if } s_i \leq p < q \leq t_i \mbox{ for all}\ 1 \leq i \leq k\\
  2 & \mbox{ otherwise }\end{array}\right.$$

It is easy to check $\phi$ is a locally complete 2-edge-colouring of $G$. 
The colouring is not unique when $cc(G) = 2$ and $G$ is of type 2 (or 3) as shown in 
Lemma \ref{lem:cc=2}. When $cc(G) = 2$ and $G$ is of type 4,
Lemma \ref{lem:cc=2} shows that the colouring $\phi$ is unique up to switching the 
colours. For all other cases, $cc(G) \geq 3$, which implies that $\overline{G}$ 
contains a triangle and hence is not bipartite. By Theorem \ref{all-in-one}, 
the colouring $\phi$ is unique up to switching the colours.
\qed 

\begin{remark} \label{rk} 
We summarize some of the facts on the edge-colourings $\phi$ of $G$ defined in 
Lemma \ref{lem:PIGcc>=3} which will be used in the rest of the section.
\begin{itemize}
\item When $cc(G) \geq 3$, $\overline{G}$ is not bipartite and hence, by 
Theorem \ref{all-in-one}, $\phi$ is unique up to switching of the colours. 
Moreover, all edges incident with $v_{s_1}$ have the same colour if and only if 
$t_1 \neq s_2$ and all edges incident with $v_{t_k}$ have the same colour 
if and only if $s_k \neq t_{k-1}$. In particular, if $G$ is of type 4 then all
edges incident with either of $v_{s_1}$ and $v_{t_k}$ are coloured with 
the same colour. 

\item When $cc(G) = 2$ and $G$ is of type 4, $G$ has no universal vertex. 
By Theorem \ref{all-in-one}, $\phi$ is unique up to switching of the colours. 
All edges incident with either of $v_{s_1}$ and $v_{t_k}$ are coloured with 
the same colour. 

\item When $cc(G) = 2$ and $G$ is of type 1, $G$ has a universal vertex and 
a locally complete 2-edge-colouring of $G$ is not unique. But according to 
Lemma \ref{lem:cc=2}, for any locally complete 2-edge-colouring of $G$, exactly 
one of $v_{s_1}$ and $v_{t_k}$ is incident with edges of the same colour and 
it is always possible to specify which one of the two vertices for having this 
property with a suitable locally complete 2-edge-colouring of $G$. 

\item For  any locally complete 2-edge-colouring of $G$ with 
$cc(G) \geq 2$ there exist edges $v_av_b, v_cv_d$ with $a \leq s_1+1$ ($=2$), 
$b > t_1$, $c < s_k$ and $d \geq t_k-1$ ($= n-1$) which are coloured with 
the same colour.
\qed
\end{itemize}
\end{remark}

The following theorem gives a full description of connected reduced proper interval 
graphs which are locally complete 2-edge-colourable.

\begin{tm} \label{thm:PIGchar}
Let $G$ be a connected reduced proper interval graph with straight ordering 
$v_1,\ldots{},v_n$. Then $G$ has has a locally complete 2-edge-colouring if and 
only if one of the following holds:
\begin{itemize}
\item[(1)] $G$ has no cutvertex, in which case either $cc(G) = 1$ or 
      $G$ is of one of the 4 types.
\item[(2)] $v_{j_1}, v_{j_2}, \dots, v_{j_c}$ are the cutvertices of $G$ where
      $1 < j_1 < j_2 < \cdots < j_c < n$, 
\begin{itemize}
\item[$\bullet$] $cc(G[1,j_1]) = 1$, or $cc(G[1,j_1]) = 2$ and $G[1,j_1]$ is of
           type 1, or $G[1,j_1]$ is of type 2 or 4,
\item[$\bullet$] $cc(G[j_c,n]) = 1$, or $cc(G[j_c,n]) = 2$ and $G[j_c,n]$ is of 
           type 1, or $G[j_c,n]$ is of type 3 or 4, and 
\item[$\bullet$] for each $i = 1, 2, \dots, c-1$, either $cc(G[j_i,j_{i+1}]) = 1$ 
           or $G[j_i,j_{i+1}]$ is of type 4.
\end{itemize}
\end{itemize}
\end{tm}
\pf Suppose that $G$ has no cutvertex. If $cc(G) = 1$ then $G = K_1$
is locally complete 2-edge-colourable. If $G$ is of one of the 4 types, then by 
Lemma~\ref{lem:PIGcc>=3} $G$ has a locally complete 2-edge-colouring.
Suppose that $G$ has cutvertices and satisfies the condition $(2)$.
By Remark \ref{rk}, $G[1,j_1]$ has a locally complete 2-edge-colouring in which
all edges incident with $v_{j_1}$ are coloured 2 and for each $i = 1, 2, \dots, c-1$,
$G[j_i,j_{i+1}]$ has a locally complete 2-edge-colouring in which all edges 
incident with either of $v_{j_i}$ and $v_{j_{i+1}}$ are coloured 1 if $i$ is odd
and 2 if $i$ is even. Finally, also by Remark \ref{rk}, $G[j_c,n]$ has a locally
complete 2-edge-colouring in which all edges incident with $v_{j_c}$ are coloured 
1 if $c$ is odd and 2 if $c$ is even. These colourings together form a locally 
complete 2-edge-colouring of $G$. 

Conversely, suppose that $G$ has a locally complete 2-edge-colouring. Then by
Proposition \ref{prop:Fhasnocol} none of $F_1, F_2, F_3$ is an induced subgraph of 
$G$. When $G$ has no cutvertex, either $cc(G) = 1$ or $G$ is of one of the 4 types
by Corollary \ref{types} $G$ is of one of the 4 types. Hence $(1)$ holds.
So assume that $G$ has cutvertices $v_{j_1}, v_{j_2}, \dots, v_{j_c}$ where 
$1 = j_0 < j_1 < j_2 < \cdots < j_c < j_{c+1} = n$. By Corollary \ref{types}, 
each $G[j_i,j_{i+1}]$ is either a clique or of one of the 4 types. 
For each $i = 1, 2, \dots, c$, since $v_{j_i}$ is a cutvertex of $G$, 
all edges of $G[j_{i-1},j_i]$ incident with $v_{j_i}$ are coloured with the same 
colour and all edges of $G[j_i,j_{i+1}]$ incident with $v_{j_i}$ are coloured with
the same colour in any locally complete 2-edge-colouring of $G$. 
By Remark \ref{rk}, each $G[j_i,j_{i+1}]$ is either a clique or of one of the 
types described in $(2)$.  
\qed\\

We close this subsection by stating a lemma which will be used in the next 
subsection. The proof follows easily from the definition of the different types.

\begin{lemma}\label{lem:pathsintypes}
Let $G$ be a proper interval graph with $cc(G) \geq 3$ which is of one of the 4 
types in Definition \ref{def:types}. Then the following holds:
\begin{itemize}
\item[(1)] If $G$ is of type 1, then 
\begin{itemize}
\item $v_{s_2}v_{t_2}v_{s_3}v_{t_3} \dots v_{s_{k-1}}v_{s_k}$ is an induced odd $(v_{s_2},v_{s_k})$-path;
\item $v_{t_2}v_{s_3}v_{t_3} \dots v_{s_{k-1}}v_{s_k}v_{t_k}$ is an induced odd $(v_{t_2},v_{t_k})$-path;
\item $v_{s_2+1}v_{t_2}v_{s_3}v_{t_3} \dots v_{s_{k-1}}v_{s_k}$ is an induced odd $(v_{s_2+1},v_{s_k})$-path;
%\item $v_{s_2-1}v_{s_2}v_{t_2}v_{s_3}v_{t_3} \dots v_{s_{k-1}}v_{s_k}v_{t_k}$ is an induced odd $(v_{s_2-1},u_{t_k})$-path;
\item $v_{s_2-1}v_{s_2+1}v_{t_2}v_{s_3}v_{t_3} \dots v_{s_{k-1}}v_{s_k}v_{t_k}$ is an induced odd $(v_{s_2-1},u_{t_k})$-path; 
\item $v_{s_1}v_{s_2}v_{t_2} \dots v_{s_{k-1}}v_{s_k}$ is an induced even $(v_{s_1},v_{s_k})$-path;
\item $v_{s_1}v_{s_2}v_{t_2}v_{s_3}v_{t_3} \dots v_{s_{k-1}}v_{s_k}v_{t_k}$ an induced odd $(v_{s_1},v_{t_k})$-path;
\item $v_{t_1}v_{t_2}v_{s_3}v_{t_3} \dots v_{s_k}v_{t_k}$ is an induced even $(v_{t_1},v_{t_k})$-path.
\end{itemize}

\item[(2)] If $G$ is of type 2, then 
\begin{itemize}
\item $v_{t_2}v_{s_3}v_{t_3} \dots v_{t_{k-1}}v_{s_k}$ is an induced odd $(v_{t_2},v_{s_k})$-path;
\item $v_{t_2}v_{s_3}v_{t_3} \dots v_{t_{k-1}}v_{s_k}v_{t_k}$ is an induced even $(v_{t_2},v_{t_k})$-path;
\item $v_{s_2+1}v_{t_2}v_{s_3}v_{t_3} \dots v_{t_{k-1}}v_{s_k}$ is an induced even $(v_{s_2+1},s_k)$-path;
\item $v_{s_1}v_{s_2}v_{t_2} \dots v_{s_{k-1}}v_{t_{k-1}}v_{s_k}$ is an induced odd $(v_{s_1},v_{s_k})$-path;
\item $v_{s_1}v_{s_2}v_{t_2}v_{s_3}v_{t_3} \dots v_{s_k}v_{t_k}$ is an induced even $(v_{s_1},v_{t_k})$-path;
\item $v_{s_2+1}v_{t_2}v_{s_3}v_{t_3} \dots, v_{s_k}v_{t_k}$ is an induced odd $(v_{s_2+1},v_{t_k})$-path.
\end{itemize}

\item[(3)] If $G$ is of type 3, then 
\begin{itemize}
\item $v_{s_2}v_{t_2}v_{s_3}v_{t_3} \dots v_{s_{k-1}}v_{s_k}$ is an induced odd $(v_{s_2},v_{s_k})$-path;
\item $v_{s_2}v_{t_2}v_{s_3}v_{t_3} \dots v_{s_{k-1}}v_{s_k}v_{s_k+1}$ is an induced even $(v_{s_2},v_{s_k+1})$-path;
\item $v_{s_1}v_{t_1}v_{s_2},v_{t_2},\dots, v_{s_k}$ is an induced odd $(v_{s_1},v_{s_k})$-path;
\item $v_{s_1}v_{t_1}v_{s_2}v_{t_2} \dots v_{s_{k-1}}v_{s_k}v_{t_k}$ is an induced even $(v_{s_1},v_{t_k})$-path;
  \item $v_{s_1}v_{t_1}v_{s_2}v_{t_2} \dots v_{s_{k-1}}$ is an induced even $(v_{s_1},v_{s_{k-1}})$-path;
\item $v_{t_1}v_{s_2}v_{t_2} \dots v_{s_{k-1}}v_{s_k-1}v_{s_k+1}v_{t_k}$ is an induced even $(v_{t_1},v_{t_k})$-path.
\end{itemize}

\item[(4)] If $G$ is of type 4, then  
\begin{itemize}
\item $v_{s_2}v_{t_2}v_{s_3}v_{t_3} \dots v_{t_{k-1}}v_{s_k}$ is an induced even $(v_{s_2},v_{s_k})$-path;
\item $v_{s_1}v_{t_1}v_{s_2}v_{t_2} \dots v_{t_{k-1}}v_{s_k}$ is an induced even $(v_{s_1},v_{s_k})$-path;
\item $v_{s_1}v_{t_1}v_{s_2}v_{t_2} \dots v_{s_k}v_{t_k}$ is an induced odd $(v_{s_1},v_{t_k})$-path;
\item $v_{t_1}v_{s_2}v_{t_2} \dots v_{t_{k-1}}v_{s_k}v_{t_k}$ is an induced even $(v_{t_1},v_{t_k})$-path.
\end{itemize}
\end{itemize}
\end{lemma}

\subsection{The case when $G$ is a proper circular-arc graph but not a proper interval graph}

In this subsection we characterize those  circular-arc graphs  which have locally complete 
2-edge-colourings. We may again restrict our study to reduced graphs. 
In view of Proposition \ref{cc2} and Theorem \ref{thm:PIGchar} we may further 
restrict it to graphs which have clique covering number at least 3 and are not 
proper interval graphs. Hence, unless otherwise mentioned, all graphs $G$ in this
subsection have clique covering number at least 3 and are reduced proper 
circular-arc graphs but not proper interval graphs.

In the proofs below we will use, without mentioning it each time, the following structural property of round orderings of proper circular-arc graphs which follows from Remark \ref{pig/pcastructure}

  \begin{lemma}\label{pcacliques}
    If $w_1,w_2,\ldots{},w_n$ is a round ordering of a proper circular-arc graph $X$, then for every edge $w_iw_j$ either $X[i,j]$ is a clique or $X[j,i]$ is a clique.
\qed
\end{lemma}

Let $v_1, v_2, \dots, v_n$ be a round ordering of $G$. Since $G$ is not a proper 
interval graph, we have $v_i \sim v_{i+1}$ for all $i$. By Lemma \ref{pcacliques}, 
if $v_i\sim v_j$, then one of $G[i,j]$ and $G[j,i]$ is complete.  

A vertex $v_i$ of $G$ is called a {\bf pseudo-cutvertex} if $v_{i-1}\nad v_{i+1}$. We call 
$G[i,j]$ a {\bf weak block} of $G$ if $v_i$ and $v_j$ are the only pseudo-cutvertices
of $G$ contained in $G[i,j]$. Each weak block of $G$ is a connected proper 
interval graph. 

We first consider the case when $G$ contains pseudo-cutvertices.

\begin{tm} \label{lem:weakblocks}
Suppose that $G$ has at least two pseudo-cutvertices. Then $G$ has a locally complete
2-edge-colouring if and only if the following hold:
\begin{itemize}
\item The number of pseudo-cutvertices is even.
\item Each weak block is either a clique or of type 4.
\end{itemize}
\end{tm}
\pf Let $v_{j_1}, v_{j_2}, \dots, v_{j_p}$ where $1\leq j_1<j_2< \cdots <j_p \leq n$
be the pseudo-cutvertices of $G$. If any of the weak blocks $G[j_i,j_{i+1}]$ 
contains one of $F_1, F_2, F_3$ as an induced subgraph, then $G$ is not locally 
complete 2-edge-colourable. Thus, {by Corollary \ref{types},}  each $G[j_i,j_{i+1}]$ is either a clique or 
of one of types 1-4. By the definition of pseudo-cutvertices and properties of 
a round 
ordering, every edge between $V(G[j_i,j_{i+1}])$ and the remaining vertices is incident 
with either $v_{j_i}$ or $v_{j_{i+1}}$. This means that any locally complete 
2-edge-colouring of $G[j_{i},j_{i+1}]$ must colour all edges of $G[j_{i},j_{i+1}]$ 
incident with $v_{j_i}$ by the same colour and all edges of $G[j_{i},j_{i+1}]$ 
incident with $v_{j_{i+1}}$ by the same colour. This can be done only when either 
$G[j_{i},j_{i+1}]$ is a clique or of type 4. Since $v_{j_{i+1}}$ is 
a pseudo-cutvertex, the colour used for the edges of $G[j_{i},j_{i+1}]$ incident 
with $v_{j_{i+1}}$ can not be the same as the colour used for the edges of 
$G[j_{i+1},j_{i+2}]$ incident with $v_{j_{i+1}}$. This means that $p$ must be 
an even number.  

Conversely, suppose that $p$ is even and each $G[j_i,j_{i+1}]$ is a clique or of 
type 4. Then by Lemma \ref{lem:PIGcc>=3} $G[j_i,j_{i+1}]$ has a locally complete 
2-edge-colouring in which all edges incident with $v_{j_i}$ or $v_{j_{i+1}}$ are 
coloured 1 when $i$ is odd and 2 when $i$ is even. These colourings together give 
a locally complete 2-edge-colouring of $G$.  
\qed
    
\begin{tm} \label{lem:onecutv}
Suppose that $G$ has exactly one pseudo-cutvertex. Then $G$ is not locally complete 
2-edge-colourable.
\end{tm}
\pf Suppose to the contrary that $G$ has a locally complete 2-edge-colouring 
$\phi: E(G) \to \{1,2\}$. Then $G$ contains none of $F_1, F_2, F_3$ and no 
odd cycle of length $\geq 5$ as an induced subgraph.

By renumbering if necessary we can assume that $v_n$ is the pseudo-cutvertex. 
Let $\ell(n) = q$ and $\gamma(n) = p$. Since $cc(G) \geq 3$, $p < q-1$.
Consider $G[1,n-1]$. It is a connected proper interval graph with straight ordering 
$v_1, v_2, \dots, v_{n-1}$ and without cutvertices. Since $v_1 \not\sim v_{n-1}$, 
$cc(G[1,n-1]) \geq 2$. Moreover, {$G[1,n-1]$} contains none of $F_1, F_2, F_3$ as 
an induced subgraph. However, $G[1,n-1]$ may contain true twins. Any true twins 
$u, v$ must have $N_G[u] = N_G[v] \cup \{v_n\}$. Such a pair can only be 
$v_p, v_{p+1}$ or $v_q, v_{q-1}$. Let $H$ be the graph obtained from $G[1,n-1]$ by 
deleting $v_{p+1}$ if $v_p, v_{p+1}$ are true twins in $G[1,n-1]$ and deleting 
$v_{q-1}$ if $v_q, v_{q-1}$ are true twins in $G[1,n-1]$. Then $H$ is a connected 
reduced proper interval graph with straight ordering 
$v_1, v_2, \dots, v_p, \dots, v_q, v_{q+1}, \dots, v_{n-1}$.  
The restriction of $\phi$ to $H$ is a locally complete 2-edge-colouring of 
$H$. 

Suppose  first that $H$ does not contain a cutvertex. By Remark \ref{rk}, there exist edges 
$v_av_b$ and $v_cv_d$ with $1 \leq a \leq 2$, $p < b < q$, $n-2 \leq c \leq n-1$, 
and $p < d < q$ such that $\phi(v_av_b) = \phi(v_cv_d)$. Since $v_nv_av_b$ and 
$v_nv_cv_d$ are induced paths in $G$, $\phi(v_nv_a) = \phi(v_nv_c)$ which implies 
$v_a \sim v_c$. We must have $v_1 \nad v_{n-2}$ as otherwise $G[1,n-2]$ and 
$G[n-1,n]$ are cliques covering $G$, a contradiction to the assumption 
$cc(G) \geq 3$. Similarly, $v_2 \nad v_{n-1}$. Hence $a = 2$ and $c = n-2$. 
So $v_2 \sim v_{n-2}$ which {together with $v_n\nad v_1$ implies that} $G[2,n-2]$ is a clique.
If $v_p \sim v_{n-1}$, then $G[n,p-1]$ and $G[p,n-1]$ are cliques covering $G$,
a contradiction to the assumption $cc(G) \geq 3$. So $v_p \nad v_{n-1}$. 
Hence $v_p = v_2$ as otherwise $v_p$ and $v_2$ are true twins in $G$, a contradiction
to the assumption that $G$ is reduced. A similar argument shows that $v_q = v_{n-2}$.
Since neither $v_2$ nor $v_{n-2}$ is a pseudo-cutvertex of $G$, $v_1 \sim v_3$ and
$v_{n-3} \sim v_{n-1}$. Since $cc(G) \geq 3$, $v_1 \nad v_{n-3}$ and 
$v_3 \nad v_{n-1}$. Hence $v_3 \neq v_{n-3}$ and $v_1v_3v_{n-3}v_{n-1}v_nv_1$
is an induced $C_5$ in $G$, a contradiction. 

Suppose now that $H$ contains cutvertices. The only possible cutvertices in $H$ 
are $v_p$ and $v_q$. This is due to the deletions of $v_{p+1}$ and $v_{q-1}$. 
We claim that if $v_p$ is a cutvertex of $H$ then $p = 2$. Indeed, suppose that
$v_p$ is a cut-vertex of $H$. Clearly, $p \geq 2$ {and we have $v_1\sim v_{p+1}$ as $v_p,v_{p+1}$ are true twins of $G[1,n-1]$}. Since $v_p$ is a cutvertex 
of $H$, $v_{p+1} \sim v_{p-1} \nad v_{p+2}$. If $p > 2$, then $v_1$ and $v_{p-1}$
are true twins in $G$, a contradiction to the assumption that $G$ is reduced.
A similar argument shows that if $v_q$ is a cutvertex of $H$ then $q = n-2$.

Assume that exactly one of $v_2, v_{n-2}$ is a cutvertex of $H$. By reversing 
the round ordering of $G$ if necessary we may assume that it is $v_2$. Then $v_1$ 
is adjacent to $v_2$ only in $H$ and $H[2,n-1]$ is a connected reduced proper 
interval graph with $cc(H[2,n-1]) \geq 2$. By Corollary \ref{types}, $H[2,n-1]$ is 
of one of the 4 types.  Since $v_2$ is a cutvertex of $H$, the restriction of $\phi$
to $H$ must colour all edges of $H[2,n-1]$ incident with $v_2$ with the same colour.
By Remark \ref{rk}
{and Lemma \ref{lem:PIGcc>=3}}, $H[2,n-1]$ is of type 3 or 4. Assume without loss of generality
that $\phi(v_2v_4) = \phi(v_2v_{\gamma(2)}) = 1$. Since every two consecutive edges 
in $v_{n-1}v_nv_2v_4$ form an induced path in $G$, $\phi(v_2v_4) = 1$ implies 
$\phi(v_{n-1}v_n) = 1$. On the other hand, $v_2v_{\gamma(2)}v_{\gamma(2)+1}$ is 
an induced path in $G$, so $\phi(v_2v_{\gamma(2)}) = 1$ implies 
$\phi(v_{\gamma(2)}v_{\gamma(2)+1})=2$. Since every two consecutive edges in
$v_{n-1}v_nv_1v_3v_{\gamma(2)}v_{\gamma(2)+1}$ form an induced path in $G$, 
$\phi(v_{\gamma(2)}v_{\gamma(2)+1}) = 2$ implies $\phi(v_{n-1}v_n) = 2$. 
This is a contradiction. The proof for the case when both $v_2$ and $v_{n-2}$ are
cutvertices of $H$ is the same. In this case we use the fact that $H[2,n-2]$
(instead of $H[2,n-1]$) is a connected reduced proper interval graph with 
$cc(H[2,n-2]) \geq 2$. {To see this, first observe that $n\geq 6$ as $G$ is not a proper interval graph. Thus $v_3$ and $v_{n-3}$ are distinct vertices. If they are adjacent, then $v_n,v_1,v_3,v_{n-3},v_{n-1}$ is an induced 5-cycle in $G$, contraction. Hence $v_3\nad v_{n-3}$, implying that $cc(H[2,n-2]) \geq 2$. }
{This completes the proof that} $G$ is not locally complete 2-edge-colourable.    
\qed\\

It remains to consider the case when $G$ has no pseudo-cutvertex.

\begin{figure}[ht]
\begin{center}
\begin{tikzpicture}[>=latex]
%\begin{pgfonlayer}{nodelayer}
\draw(0,0) circle(3);

\node [label={left:$v_1$}] [style=blackvertex] (1) at (-3,0) {};
\node [label={above:$v_a$}] [style=blackvertex] (2) at (-1,2.8) {};
\node [label={below:$v_b$}] [style=blackvertex] (3) at (-1,-2.8) {};

\node [label={left:$v_2$}] [style=blackvertex] (4) at (-2.9,.9) {};
\node [label={left:$v_n$}] [style=blackvertex] (5) at (-2.9,-.9) {};

\node [label={left:$v_{a-1}$}] [style=blackvertex] (7) at (-1.8,2.4) {};
\node [label={left:$v_{b+1}$}] [style=blackvertex] (8) at (-1.8,-2.4) {};

\draw[dashed] (1) -- (2) -- (3) -- (1);

\draw[-] (1) to (7) {};
\draw[-] (1) to (8) {};

\draw[-] (2) to (4) {};
\draw[-] (3) to (5) {};

\end{tikzpicture}
\end{center}
\vspace{-2mm}
\caption{Illustrating the choice of vertices $v_1, v_a, v_b$ in $G$.}
\label{1ab}
\end{figure}
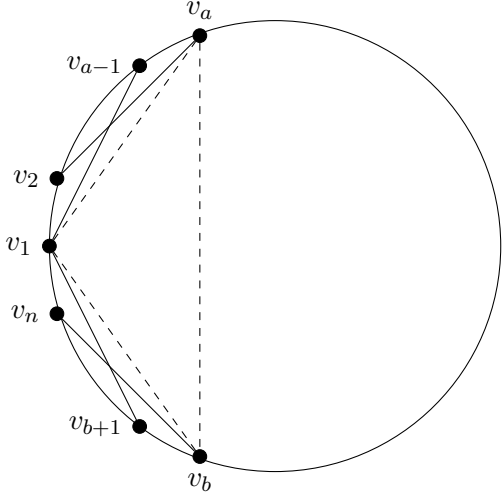

%\begin{prp}\label{lem:nopseudocv}
\begin{tm}\label{lem:nopseudocv}
Suppose that $G$ has no pesudo-cutvertex and is perfect. Then $G$ has is locally 
complete 2-edge-colourable if and only if it does not contain any of $F_1, F_2, F_3$
as an induced subgraph.
%\end{prp}
\end{tm}
\pf The necessity follows from Proposition \ref{prop:Fhasnocol}. For sufficiency,
suppose that $G$ is a reduced graph which is perfect and does not contain any of 
$F_1, F_2, F_3$ as 
an induced subgraph. In all of the proof below we let $v_1,v_2,\ldots{},v_n$ be a round ordering of $G$.
As $cc(G)\geq 3$, $\overline{G}$ is not bipartite so it 
contains an odd cycle. The shortest odd cycle is induced and since $G$ is perfect, 
$\overline{G}$ contains no induced odd cycle of length at least 5. Hence 
$\overline{G}$ has a 3-cycle, that is,
$G$ contains 3 pairwise non-adjacent vertices. 

Without loss of generality assume that the  vertices $v_1, v_a, v_b$ are pairwise non-adjacent chosen in such a way that $G[b,a]$ has as few vertices as possible.

The following is a direct consequence of Lemma \ref{pcacliques}.
\begin{eqnarray}
  %& v_1\mbox{ is adjacent to every vertex in $G[2,a-1]$ and in $G[b+1,n]$.}\\
  & v_1\mbox{ is not adjacent to any  vertex in $G[a,b]$.\label{v1nad}}\\
  & v_a\mbox{ is not adjacent to any vertex in $G[b,1]$.\label{vanad}}\\
  & v_b\mbox{ is not adjacent to any vertex in $G[1,a]$.\label{vbnad}}
%\\\nonumber
  \end{eqnarray}

We start with  a very useful property of this way of choosing the vertices $v_1,v_a,v_b$. we shall use this property many times in the proof of the theorem.
See Figure~\ref{1ab} for illustration.

\begin{claim}
  \label{vav2vbvn}
  The following holds for the vertices $v_1,v_a,v_b$
  \begin{eqnarray}
    v_1\sim v_{b+1}\\
    v_1\sim v_{a-1}\\
    v_a\sim v_2 &\\
    v_n\sim v_b
    \end{eqnarray}
\end{claim}
\pf
  As $v_a\nad v_1$ and $v_a\nad v_b$ we see from (\ref{vanad}) that $v_a\nad v_j$ for every vertex $v_j$ in $G[b,n]$.
  In particular $v_a\nad v_{b+1}$ so as $G[b+1,a]$ has fewer vertices than $G[b,a]$ we must have $v_1\sim v_{b+1}$. By a similar argument (using that $v_b\nad v_i$ for every vertex $v_i$ in $G[1,a]$) we see that $v_1\sim v_{a-1}$. Suppose $v_a\nad v_2$. Then we must have $v_2\sim v_{b+1}$ or 
 $v_1,v_a,v_{b+1}$ would be  a better choice of pairwise non-adjacent vertices. Now we see that $v_1$ and $v_2$ are true twins of $G$, contradicting the assumption that $G$ is reduced. Hence $v_a\sim v_2$. Suppose $v_n\nad v_b$. then $v_n\nad v_{a-1}$ as $v_1,v_n$ are not true twins of $G$, but now $v_b,v_n,v_{a-1}$ pairwise non-adjacent and form a  better choice than $v_1,v_a,v_b$, contradiction. Hence we have $v_n\sim v_b$.
 \qed\\

As a consequence of Claim \ref{vav2vbvn} we have
\begin{equation}\label{ell2gammab}
  \ell(a)=2 \mbox{ and } \gamma(b)=n
  \end{equation}

\vspace{2mm}

Consider $G[1,b]$.  It  is a connected proper interval graph with straight 
ordering $v_1, v_2, \dots, v_b$. Let $H$ be the reduced subgraph of $G[1,b]$. 
Since $G$ does not contain any of $F_1, F_2, F_3$ as an induced subgraph, $H$ does 
not contain any of $F_1, F_2, F_3$ as an induced subgraph. 
Assume that $v_{j_1}, v_{j_2}, \dots, v_{j_q}$ where 
$1 \leq j_1 < j_2 < \cdots < j_q \leq b$ are the vertices of $H$. 
Since $v_1, v_a, v_b$ are pairwise non-adjacent in $G$, they are not true twins of each
other so we can assume that $H$ contains all of them. This implies 
$cc(H) \geq 3$, $j_1 = 1$, $j_q = b$, and $j_p = a$ for some $p$. 
By \ref{vbnad}, $v_{a-1} \nad v_b$. Thus, as we have $v_1 \sim v_{a-1} \sim v_a \nad v_1$, the vertex $v_{a-1}$ is not
a twin with any of $v_1, v_a, v_b$ in $G[1,b]$. So we can also assume that 
$v_{a-1}$ is a vertex in $H$ which means that $j_{p-1} = a-1$.
To simplify subscripts, we rename the vertices of $H$ by letting 
$u_i = v_{j_i}$ for all $1 \leq i \leq q$. With the new names,
$H[i,j]$ is the subgraph of $H$ induced by vertices $u_i, u_{i+1}, \dots, u_j$.\\

{\noindent{}{\bf Case 1: $H$ has no cutvertex}}

%Suppose first that $H$ contains no cutvertex.
By Corollary \ref{types}, $H$ is 
one of the 4 types. Let $H[s_i,t_i]$ ($1 \leq i \leq k$) be the canonical clique 
covering of $H$. Then $u_{s_1} = v_1$, $u_{t_1} = v_{a-1}$ and $u_{t_k} = v_b$.
Since $H$ is an induced subgraph of $G$ and we only removed a vertex if it was a twin of some remaining vertex, we can assume without loss of generality that for all $i\in [k-1]$ the vertices $u_{t_i}$ and $u_{s_{i+1}}$ are
either the same vertex of $G$ or consecutive vertices in the round ordering of $G$.\\

We shall use the following properties many times without mentioning it every time. 

\begin{claim}
  \label{nonadj}
When there are at least 3 cliques in the canonical clique cover of $H[1,b]$ (that is, $k\geq 3$) the following holds.
\begin{itemize}
\item[(a)] No vertex of $G[b,1]$ is adjacent to a vertex of $G[v_a,u_{s_k-1}]$
  \item[(b)] No vertex of $G[1,a]$ is adjacent to a vertex of $G[u_{t_2+1},v_b]$
\end{itemize}
\end{claim}
\pf
  To prove (a) suppose $v_j\sim v_c$ for some $v_j$ in $G[v_b,v_1]$ and  $v_c$ in $G[v_a,u_{s_k-1}]$. As $v_a\nad v_1$ it follows from 
  Lemma \ref{pcacliques} that $G[v_c,v_j]$ is a clique. However this implies that $\ell(u_{t_k})$ is a vertex before $u_{s_k}$ in the round ordering of $V(G)$, contradicting  Lemma \ref{canonical} (2). Hence (a) holds.
  To prove (b) suppose that $v_i\sim v_d$ for some $v_i$ in $G[v_1,v_a]$ and $v_d$ in $G[u_{s_k+1},v_b]$. As $v_1\nad v_b$ it follows from Lemma \ref{pcacliques} that $G[v_i,v_d]$ is a clique and hence $\gamma(u_{s_2})$ is a vertex after $u_{t_2}$ in the round ordering of $V(G)$. Again this contradicts Lemma \ref{canonical} (2). Hence (b) holds.
  \qed

\medskip

\begin{figure}[hbt]
\begin{center}
\begin{tikzpicture}[>=latex]
%\begin{pgfonlayer}{nodelayer}
\draw(0,0) circle(3);

\node [label={left:$v_1$}] [style=blackvertex] (1) at (-3,0) {};
\node [label={above:$v_a$}] [style=blackvertex] (2) at (-1,2.8) {}; 
\node [label={below:$v_b$}] [style=blackvertex] (3) at (-1,-2.8) {};

\node [label={right:$u_{t_3}$}] [style=blackvertex] (4) at (2.8,1.2) {};
\node [label={above:$u_{t_2}$}] [style=blackvertex] (5) at (1,2.8) {};
\node [label={below:$u_{s_k}$}] [style=blackvertex] (6) at (1,-2.8) {};

\node [label={left:$v_{a-1}$}] [style=blackvertex] (7) at (-1.8,2.4) {};
\node [label={left:$v_{b+1}$}] [style=blackvertex] (8) at (-1.8,-2.4) {};

\node [label={right:$u_{s_{k-1}}$}] [style=blackvertex] (9) at (2.5,-1.7) {};
\node [label={right:$u_{s_3}$}] [style=blackvertex] (10) at (1.8,2.4) {};
  
\draw[dashed] (1) -- (2) -- (3) -- (1) {};

\draw[-] (1) to (8) {};
\draw[-] (1) to (7) {};

\draw [ ] (7) to [out=-80,in=-80] (5) {};
\draw [ ] (10) to [out=240,in=240] (4) {};
\draw [ ] (3) to [out=90,in=90] (6) {}; 
\draw [ ] (6) to [out=90,in=100] (9) {}; 

\node [style=textbox]  at (0,-3.5) {$(I)$};

\draw(8,0) circle(3);

\node [label={left:$v_1$}] [style=blackvertex] (1) at (5,0) {};
\node [label={above:$v_a$}] [style=blackvertex] (2) at (7,2.8) {}; 
\node [label={below:$v_b$}] [style=blackvertex] (3) at (7,-2.8) {};

\node [label={right:$u_{t_3}$}] [style=blackvertex] (4) at (10.8,1.2) {};
\node [label={above:$u_{t_2}$}] [style=blackvertex] (5) at (9,2.8) {};
\node [label={below:$u_{s_k}$}] [style=blackvertex] (6) at (9,-2.8) {};

\node [label={right:$u_{t_{k-1}}$}] [style=blackvertex] (11) at (9.5,-2.6) {};

\node [label={left:$v_{a-1}$}] [style=blackvertex] (7) at (6.2,2.4) {};
\node [label={left:$v_{b+1}$}] [style=blackvertex] (8) at (6.2,-2.4) {};

\node [label={right:$u_{s_{k-1}}$}] [style=blackvertex] (9) at (10.5,-1.7) {};
\node [label={right:$u_{s_3}$}] [style=blackvertex] (10) at (9.8,2.4) {};

\draw[dashed] (1) -- (2) -- (3) -- (1) {};
                        
\draw[-] (1) to (8) {}; 
\draw[-] (1) to (7) {};
                        
\draw [ ] (7) to [out=-80,in=-80] (5) {};
\draw [ ] (10) to [out=240,in=240] (4) {};
\draw [ ] (3) to [out=90,in=90] (6) {};
\draw [ ] (11) to [out=90,in=100] (9) {};

\node [style=textbox]  at (8,-3.5) {$(II)$};

\draw(0,-8) circle(3);

\node [label={left:$v_1$}] [style=blackvertex] (1) at (-3,-8) {};
\node [label={above:$v_a$}] [style=blackvertex] (2) at (-1,-5.2) {}; 
\node [label={below:$v_b$}] [style=blackvertex] (3) at (-1,-10.8) {};

\node [label={right:$u_{t_3}$}] [style=blackvertex] (4) at (2.8,-6.8) {};
\node [label={above:$u_{t_2}$}] [style=blackvertex] (5) at (1,-5.2) {};
\node [label={below:$u_{s_k}$}] [style=blackvertex] (6) at (1,-10.8) {};

\node [label={left:$v_{a-1}$}] [style=blackvertex] (7) at (-1.8,-5.6) {};
\node [label={left:$v_{b+1}$}] [style=blackvertex] (8) at (-1.8,-10.4) {};

\node [label={right:$u_{s_{k-1}}$}] [style=blackvertex] (9) at (2.5,-9.7) {};
\node [label={right:$u_{s_3}$}] [style=blackvertex] (10) at (1.8,-5.6) {};

\draw[dashed] (1) -- (2) -- (3) -- (1) {};
                        
\draw[-] (1) to (8) {}; 
\draw[-] (1) to (7) {};
                        
\draw [ ] (2) to [out=-80,in=-80] (5) {};
\draw [ ] (10) to [out=240,in=240] (4) {};
\draw [ ] (3) to [out=90,in=90] (6) {};
\draw [ ] (6) to [out=90,in=100] (9) {};

\node [style=textbox]  at (0,-11.5) {$(III)$};

\draw(8,-8) circle(3);

\node [label={left:$v_1$}] [style=blackvertex] (1) at (5,-8) {};
\node [label={above:$v_a$}] [style=blackvertex] (2) at (7,-5.2) {};
\node [label={below:$v_b$}] [style=blackvertex] (3) at (7,-10.8) {};

\node [label={right:$u_{t_3}$}] [style=blackvertex] (4) at (10.8,-6.8) {};
\node [label={above:$u_{t_2}$}] [style=blackvertex] (5) at (9,-5.2) {};
\node [label={below:$u_{s_k}$}] [style=blackvertex] (6) at (9,-10.8) {};

\node [label={left:$v_{a-1}$}] [style=blackvertex] (7) at (6.2,-5.6) {};
\node [label={left:$v_{b+1}$}] [style=blackvertex] (8) at (6.2,-10.4) {};

\node [label={right:$u_{s_{k-1}}$}] [style=blackvertex] (9) at (10.5,-9.7) {};
\node [label={right:$u_{s_3}$}] [style=blackvertex] (10) at (9.8,-5.6) {};

\node [label={right:$u_{t_{k-1}}$}] [style=blackvertex] (11) at (9.5,-10.6) {};

\draw[dashed] (1) -- (2) -- (3) -- (1) {};

\draw[-] (1) to (8) {};
\draw[-] (1) to (7) {};

\draw [ ] (2) to [out=-80,in=-80] (5) {};
\draw [ ] (10) to [out=240,in=240] (4) {};
\draw [ ] (3) to [out=90,in=90] (6) {};
\draw [ ] (11) to [out=90,in=100] (9) {};

\node [style=textbox]  at (8,-11.5) {$(IV)$};

\end{tikzpicture}
\end{center}
\vspace{-2mm}
\caption{Illustrating the 4 types of $H$.}
\label{cases}
\end{figure}
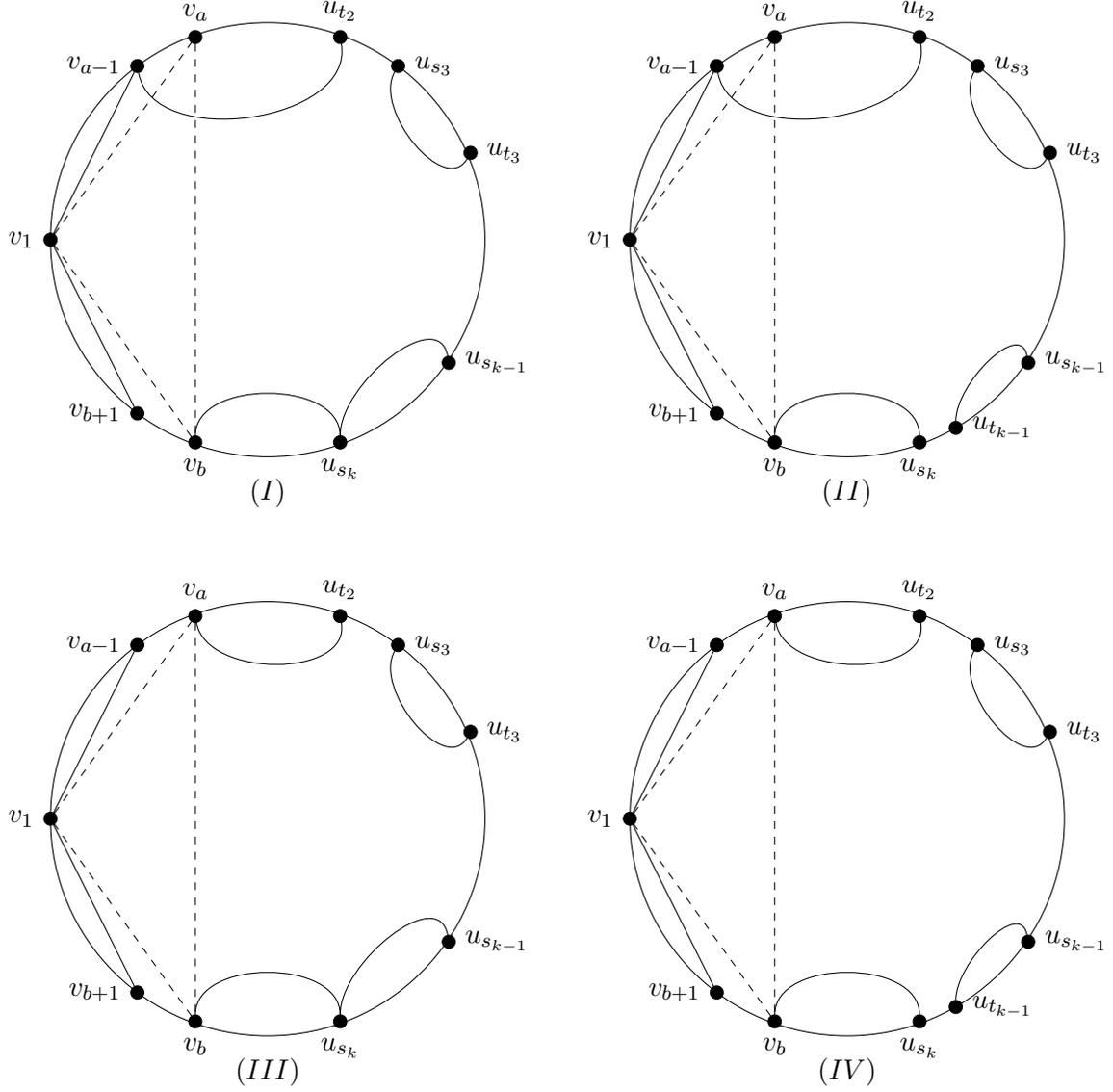

\noindent{}{\bf Case 1.1: $H$ is of type 1} (see Figure \ref{cases}$ (I)$)

By definition of type 1, we have $u_{t_1} = u_{s_2}$, $u_{t_{k-1}} = u_{s_k}$
and $u_{t_i} = u_{s_{i+1}-1}$ for all $2 \leq i \leq k-2$. 
Note that $k>2$ as $u_{t_1}=v_{a-1}\nad v_b=u_{t_k}$.

By Lemma \ref{lem:pathsintypes}, there exists an induced odd $(u_{s_2},u_{s_k})$-path 
$P_1$, an induced odd $(u_{t_2},u_{t_k})$-path $P_2$, an induced odd
$(u_{s_2+1},u_{s_k})$-path $P_3$, and an induced odd $(u_{s_2-1},u_{t_k})$-path
$P_4$ in $H$. Using these paths, we obtain an induced odd cycle in $G$ as follows:
  \begin{itemize}
\item If $v_{a-1} \sim v_{b+1} \sim u_{s_k}$ and $k>3$ then it follows from Claim \ref{nonadj} that $v_{b+1}P_1v_{b+1}$ is
  an induced odd cycle of length $\geq 5$ in $G$. {If $k=3$ then  we have 
$u_{s_2}=u_{t_1}$ and $u_{s_3}=u_{t_2}$.
    Then $v_{b+1}v_{a-2}v_au_{s_3-1}u_{s_3+1}v_{b+1}$ is an induced 5-cycle. Here we used Lemma \ref{2-canonical} (1) to see that $v_{a-2}\nad u_{s_3-1}$ and $v_a\nad u_{s_3+1}$ and the fact that $\ell(u_{t_i})=u_{s_i}$ to see that $v_{b+1}\nad u_{s_3-1}$ and $v_{a-2}\nad u_{s_3+1}$.}
\item If $v_{a-1} \nad v_{b+1} \nad u_{s_k}$, then $v_1v_{a-1}P_2v_{b+1}v_1$ 
is an induced odd cycle of length $\geq 5$ in $G$.
\item If $v_{a-1} \nad v_{b+1} \sim u_{s_k}$, then $v_1v_{a-1}P_3v_{b+1}v_1$
is an induced odd cycle of length $\geq 5$ in $G$.
\item If $v_{a-1} \sim v_{b+1} \nad u_{s_k}$, then $v_{b+1}P_4v_{b+1}$ 
is an induced odd cycle of length $\geq 5$ in $G$.
\end{itemize}
All these contradict the assumption that $G$ is perfect.\\

\noindent{}{\bf Case 1.2: $H$ is of type 2} (see Figure~\ref{cases}$ (II)$)

By definition of type 2, we have  $u_{t_1} = u_{s_2}$ and 
$u_{t_i} = u_{s_{i+1}-1}$ for all $2 \leq i \leq k-1$.
By Lemma \ref{lem:pathsintypes}, there exists an induced odd 
$(u_{t_2},u_{s_k})$-path 
$P_1$ and an induced even $(u_{t_2},u_{t_k})$-path $P_2$, and
an induced even $(u_{s_2+1},s_k)$-path $P_3$ in $H$.

\begin{claim}
  $v_{a-1} \nad v_{b+1} \nad u_{s_k}$ %and $v_{a-2} \sim v_{b+1}$.
\end{claim}
\pf
   If $v_{a-1} \nad v_{b+1} \sim u_{s_k}$, then $v_1v_{a-1}P_1v_{b+1}v_1$
is an induced odd cycle  of length $\geq 5$ in $G$. 
If $v_{a-1} \sim v_{b+1} \nad u_{s_k}$, then $v_{a-1}P_2v_{b+1}v_{a-1}$
is an induced odd cycle of length $\geq 5$ in $G$.
If $v_{a-1} \sim v_{b+1} \sim u_{s_k}$, then $v_2P_3v_{b+1}v_2$
  is an induced odd cycle of length $\geq 5$ in $G$.
  All of these contradict that $G$ is perfect so the claim holds.
 \qed

If we also have  $v_{a-2} \nad v_{b+1}$, then
$v_1v_{a-2}v_aP_2v_{b+1}v_1$ is an induced odd cycle of length $\geq 5$ in $G$.

 so we must have
$v_{a-2} \sim v_{b+1}$.

\begin{itemize}
  \item By Lemma \ref{2-canonical} (1) there is no vertex $v$ in $G[v_{a-1},u_{t_2}]$ which is adjacent to $v_{a-2}$ and 
$u_{s_3}$.% as otherwise the vertices $v_{b+1}, v_{a-2},v_{a-1}, v, u_{t_2}, u_{s_3}, u_{t_3}$
%induce a copy of $F_1$, a contradiction. 
%There is no vertex $v$ in $G[v_{b+1},v_{a-2}]$ which is adjacent to $v_b$ and 
%$v_{a-1}$, as otherwise $u_{t_2}, v_{a-1}, v_{a-2}, v, v_{b+1}, v_b, u_{s_k}$ 
%induce a copy of $F_1$, a contradiction.
\item By (\ref{vbnad}) we have $v_b\nad v_j$ for every $j\in [a-1]$
\item There is no vertex $w$ in $G[u_{s_k},u_{t_k}]$ such that $w\sim v_{b+1}$ and $w\sim u_{t_{k-1}}$ since otherwise $wv_{b+1}v_{a-2}u_{s_2}u_{t_2}u_{s_3}\ldots{}u_{t_{k-1}}w$ would be an induced odd cycle of length $\geq 5$.
  \item There is  no vertex $z$ in $G[v_{b+1},v_n]$ such that $z\sim u_{s_2}=v_{a-1}$ since otherwise 
    $zu_{s_2}u_{t_2}u_{s_3}\ldots{}u_{s_k}u_{t_k}z$ would be an induced odd cycle of length $\geq 5$. Hence $\ell(a-1)=1$.
    \end{itemize}
  
Define an edge colouring $\phi: E(G) \to \{1,2\}$ as follows: 
$$\phi{}(xy)=\left\{\begin{array}{ll}
  1 & \mbox{ if } x , y\ \mbox{are\ vertices\ in}\ G[b+1,a-2]\ \mbox{or}\ 
    G[u_{s_i},u_{t_i}]\  \mbox{for each}\ 2 \leq i \leq k\\
    2 & \mbox{ otherwise }\end{array}\right.$$

Using the observations above, one can easily check that $\phi$ is a locally complete 2-edge-colouring. \\

\noindent{}{\bf Case 1.3:  $H$ is of type 3} (see Figure~\ref{cases}$ (III)$)

  By definition of type 3, we have $u_{t_{k-1}} = u_{s_k}$ and 
$u_{t_i} = u_{s_{i+1}-1}$ for all $1 \leq i \leq k-2$. 

By Lemma \ref{lem:pathsintypes}, there exists an induced odd $(u_{s_2},u_{s_k})$-path
$P_1$ and an induced even $(u_{s_2},u_{s_k+1})$-path $P_2$.

\begin{claim}
  $v_{a-1} \nad v_{b+1} \nad u_{s_k}$
\end{claim}
\pf
If $v_{a-1} \nad v_{b+1} \sim u_{s_k}$, then $v_1v_{a-1}P_1v_{b+1}v_1$ is
an induced odd cycle of length $\geq 5$ in $G$.
If $v_{a-1} \sim v_{b+1} \nad u_{s_k}$, then $v_{a-1}P_1t_kv_{b+1}v_{a-1}$
an induced odd cycle of length $\geq 5$ in $G$. Here we used that $v_a\nad v_{b+1}$ by (\ref{vanad}) and $u_{t_k}=v_b\nad v_{a-1}$ by (\ref{vbnad}).
If $v_{a-1} \sim v_{b+1} \sim u_{s_k}$, then $v_{a-1}P_2v_{b+1}v_{a-1}$ is 
an induced odd cycle of  length $\geq 5$ in $G$.
Each of these contradict the assumption that $G$ is perfect, so we must have 
$v_{a-1} \nad v_{b+1} \nad u_{s_k}$. \\

Now we list some properties which finally lead to a contradiction.
\begin{itemize}
\item We must have 
$v_{b+1} \sim u_{s_k+1}$, as otherwise $v_1v_{a-1}P_2v_bv_{b+1}v_1$ is an induced odd
cycle of length $\geq 5$ in $G$, a contradiction. 
\item Since $v_{b+1}$ is not a pseudo-cutvertex of $G$, $n > b+1$. 
\item $v_n\sim u_{s_k+1}$ as otherwise If 
$v_1, v_n, v_{b+1}, v_b, u_{s_k+1}, u_{s_k}, u_{s_k-1}$ 
induce a copy of $F_2$ in $G$, a contradiction to the assumption.

%\item Since $v_{a-1}$ is not a pseudo-cutvertex of $G$, $1 < \ell(a)=2 < a-1$. 
\item As $v_1$ is not a pseudo-cutvertex of $G$ $v_2\sim v_n$ so now we see that $v_2P_2v_nv_2$ is
an induced odd cycle of length $\geq 5$ in $G$, a contradiction.

\end{itemize}

\noindent{}{\bf Case 1.4:  $H$ is of type 4} (see Figure~\ref{cases}$ (IV)$)

By definition of type 4, we have $u_{t_i} = u_{s_{i+1}-1}$ for all 
  $1 \leq i \leq k-1$. Note that $k>2$ must hold as $v_a\nad v_b$. This implies that $u_{t_1}=v_{a-1}\nad u_{s_k}$\\
  By Lemma \ref{lem:pathsintypes}, there exists an induced even 
  $(u_{s_2},u_{s_k})$-path $P_1$.\\
  
If $v_{a-1} \sim v_{b+1} \sim u_{s_k}$, then 
$v_{a-1}P_1v_{b+1}v_{a-1}$ is an induced odd cycle of length $\geq 5$.
If $v_{a-1} \nad v_{b+1} \nad u_{s_k}$, then 
  $v_1v_{a-1}P_1v_bv_{b+1}v_1$ is an induced odd cycle of length $\geq 5$.\\

Suppose that $v_{a-1} \sim v_{b+1} \nad u_{s_k}$. Now we make the following observations.
\begin{itemize}
  
\item There is no $v$ in 
$G[u_{s_k},u_{t_k}]$ which is adjacent to both $u_{t_{k-1}}$ and $v_{b+1}$. Suppose there is such a vertex, then if $k>3$ the  vertices $v_{a-1}, v_{b+1}, v_b, v, u_{s_k}, u_{t_{k-1}}, u_{s_{k-1}}$ induce
a copy of $F_1$, a contradiction. So $k=3$ and then $v_{a-1}u_{s_2}u_{t_2}vv_{b+1}v_{a-1}$ is an induced 5-cycle, contradiction.
\item %Since $v_1\nad v_a$ and $v_b\nad v_1$
  By (\ref{vanad}) and (\ref{vbnad}) there is no vertex  $x$ in $G[v_{b+1},v_{a-1}]$ such that $x\sim v_a$ and $x\sim v_b$.
  \item By Lemma \ref{2-canonical} (1)  there is no  vertex $w$ of $G[u_{s_2},u_{t_2}]$ such that $v_{a-1}\sim w\sim u_{s_3}$.
\end{itemize}

Define an edge colouring $\phi: E(G) \to \{1,2\}$ as follows:
\[\phi{}(xy)=\left\{\begin{array}{ll}
  1 & \mbox{ if } x , y\ \mbox{are\ vertices\ in}\ G[b+1,a-1]\ \mbox{or}\
    G[u_{s_i},u_{t_i}]\  \mbox{for\ each}\ 2 \leq i \leq k\\
  2 & \mbox{ otherwise }\end{array}\right.
\]
The observations above imply that $\phi$ is a locally complete 2-edge-colouring of 
$G$. \\

The only remaining possibility is $v_{a-1} \nad v_{b+1} \sim u_{s_k}$. 
Recall that $\ell(a)=2$ and  that  $\gamma(b)=n$ by (\ref{ell2gammab}).
Since $v_1$ is not a pseudo-cutvertex of $G$ we have $v_2\sim v_n=v_{\gamma(b)}$. We now derive a series of properties that will finally allow us to obtain a locally complete 2-edge-colouring of $G$.

\begin{itemize}
\item If $v_2\nad v_{b+1}$, as otherwise  $v_2P_1v_{b+1}v_2$ would be an odd induced odd cycle of length at least 5, contradiction.
  \item If $v_n\sim v_{a-1}$, as otherwise 
$v_a, v_{a-1}, v_2, v_1, v_n, v_{b+1}, v_b$ induce a copy of
    $F_2$, a contradiction to the assumption.
  \item If $v_n\nad u_{s_k}$, as otherwise  $v_{a-1}P_1v_nv_{a-1}$ is an odd induced cycle, contradiction.
  \item $\ell(2)=\ell(a-1)$ since otherwise $v_a,v_{a-1},v_2,v_1,v_n,v_{\ell(2)},v_b$ induce a copy of $F_3$.
  \item $u_{s_k}\nad v_{\ell(2)}$ as otherwise $v_{\ell(2)}v_{a-1}P_1v_{\ell(2)}$ would be an induced odd cycle of length at least 5.
   \item $u_{s_k}\sim v_{\ell(2)-1}$ as otherwise $v_bv_{\ell(2)-1}v_1v_{a-1}P_1v_b$ would be an induced odd cycle of length at least 5.
%\item As  $v_{a-1}$ and $v_2$ are not true twins there exists a vertex $w$ of $G[u_{s_2},u_{t_2}]$ such that $v_{l(a)}\nad w\sim v_{a-1}$.

%\item For every vertex $w$ as above we have $w\nad u_{s_3}$, since otherwise $wu_{s_3}u_{t_3}\ldots{}u_{s_k}v_{b+1}v_{\gamma(b)v_{a-1}w$ is an induced odd cycle, contradiction.
   \item By the properties of the cannonical clique cover $u_{t_k}\nad u_{t_{k-1}}$ so $v_{b+1}\nad u_{t_{k-1}}$.
   \item There is no vertex $z$ in $G[u_{s_k},u_{t_k}]$ such that $z\sim v_{\ell(2)}$ and $z\sim u_{t_k-1}$ since then
     $zv_{\ell(2)}v_{a-1}P_1[u_{s_2},u_{t_{k-1}}]z$ would be an odd induced cycle of length at least 5 in $G$,
   \item By Lemma \ref{2-canonical} there is no vertex $w\in G[u_{s_2},u_{t_2}]$ so that $v_{a-1}\sim w\sim u_{s_3}$.
     \item As $\ell(a)=2$ there is no vertex $v_f$ in $G[v_{\ell(2)},v_{a-1}]$ such that $v_f\sim v_{\ell(2)-1}$ and $v_f\sim v_a$.
    
\end{itemize}

\noindent{} Define $$\phi{}(xy)=\left\{\begin{array}{ll}
  1 & \mbox{ if } x , y\ \mbox{are\ vertices\ in}\ G[u_{s_k},v_{\ell(2)-1}]\ \mbox{or}\
  G[v_{\ell(2)},v_{a-1}]\  \mbox{ or } G[u_{s_i},u_{t_i}] \mbox{ for some  }\ 2 \leq i \leq k-1\\
  2 & \mbox{ otherwise }\end{array}\right.$$

It follows from the observations above that $\phi$ is a locally complete 2-edge-colouring of $G$.\\

\noindent{}{\bf Case 2: $H$ has cut-vertices}. 

  As $G$ is reduced and has no pseudo-cutvertices, it is easy to check that the  only possible cutvertices of $H$ are
$u_2$ (in which case $u_2 = v_{a-1}$ and $u_3 = v_a$) or $u_{q-1}$. 
When $u_{q-1}$ is a cutvertex of $H$, we can assume that $u_{q-1}$ is a vertex in
$G$ adjacent to $v_{b+1}$, because $v_b$ is not a pseudo-cutvertex of $G$. Furthermore, as $u_{q-1}$ is not a pseudo-cutvertex of $G$, its predecessor $u^-_{q-1}$ in the round ordering of $V(G)$ is adjacent to $v_b$. Note that when $u_2$ is a cutvertex of $H[1,b]$ the vertex $v_2$ must be a twin of $v_{a-1}$ in $H[1,b]$ since $v_2\sim v_a\nad v_1$.

If both $u_2$ and $u_{q-1}$ are cutvertices of $H$, then we claim that $v_{a-1}=u_2\nad u_{q-1}$.  Suppose $v_{a-1}\sim u_{q-1}$.
As $u_{q-1}$ and $u^-_{q-1}$ are not true twins of $G$, there is a vertex $z$ in $G[v_{b+1},v_n]$ such that $u^-_{q-1}\nad z\sim u_{q-1}$. This implies that $v_2\sim z$ as $v_1v_2u^-_{q-1}v_bzv_1$ is not an induced cycle of $G$. Hence  $w=\ell{}(v_2)$ is a vertex of $G[v_{b+1},z]$. Since $v_2$ and $v_{a-1}$ are true twins of $H[1,b]$ but not of $G$ we have $v_{a-1}\nad w$.  Furthermore  $v_{a-1}\sim v_n$ as otherwise $v_1v_{a-1}u^-_{q-1}v_bv_nv_1$ would be an induced 5-cycle in $G$. Now we see that the vertices $v_b,w,v_n,v_1,v_2,v_{a-1},v_a$ induce a copy of $F_3$ in $G$, contradiction. Hence we must have $u_2\nad u_{q-1}$ and thus $cc(H[2,q-1])\geq 2$.

If $u_2$ is a cutvertex of $H$ but $u_{q-1}$ is not a cutvertex of $H$,
then $cc(H[2,b])\geq 2$ since it follows from (\ref{vbnad}) that $v_{a-1}\nad v_b$. 

Finally, if $u_2$ is not a cutvertex of $H$ but $u_{q-1}$ is, then  $cc(H[1,q-1])\geq 2$ as $v_1\nad u_{q-1}$ by (\ref{v1nad}).

Hence in all three cases we can assume by Corollary \ref{types} that the respective  reduced proper interval graphs $H[2,q-1]$, $H[2,q]$, $H[1,q-1]$ are of one of 
the 4 types.\\

  \noindent{}{\bf Case 2.1: $u_2$ is a cutvertex of $H$.}

  We first study the cases when $u_2=v_{a-1}$ and possibly also $u_{q-1}$ is a cutvertex of $H$.
  As mentioned above, in that case the vertex $v_2$ must be a twin of $v_{a-1}$ in $G[1,b]$ as $v_2\sim v_a\nad v_1$. This implies that $v_2\sim u_{t_1}$  no matter which of the 4 types $H[2,q-1]$ or $H[2,q]$ has. This again implies   that $v_{a-1}\nad v_{b+1}$ must hold as $v_2$ and $v_{a-1}$ are not true twins of $G$.  Recall that we always have $\ell{}(u_{t_k})=u_{s_k}$ and when we consider $H[2,q-1]$ we have $u_{t_k}=u_{q-1}$ and when we consider $H[2,q]$ we have $u_{t_k}=u_q=v_b$. In both cases we have $u_{t_k}\sim v_{b+1}$. \\

  The following claim, which we use several times below without mentioning it, is proved similarly to the way  we proved Claim \ref{nonadj}.

  \begin{claim}
  \label{nonadj2}
When there are at least 3 cliques in the canonical clique cover of $H[2,b]$ resp. $H[1,q-1]$  (that is, $k\geq 3$) the following holds.
\begin{itemize}
\item[(a)] No vertex of $G[b,1]$ is adjacent to a vertex of $G[v_a,u_{s_k-1}]$
  \item[(b)] No vertex of $G[1,a-1]$ is adjacent to a vertex of $G[u_{t_1+1},v_b]$
\qed
\end{itemize}
  \end{claim}

  \begin{claim}   
    \label{ut1nonadj}
    $u_{t_1}\nad v_{b+1}$ and $u_{t_1}\nad v_b$.
  \end{claim}
  \pf It follows from Lemma \ref{canonical} (2) that $\ell{}(u_{t_2})=u_{s_2}$ is either $u_{t_1}$ or $u_{t_1+1}$ so we cannot have $v_{b+1}\sim u_{t_1-1}$ since that would imply that $u_{t_2}\sim u_{t_1-1}$. Now we see  that $u_{t_1}\nad v_{b+1}$ as otherwise $v_{b+1}v_1v_{a-1}u_{t_1-1}u_{t_1+1}v_{b+1}$ would be an induced 5-cycle in $G$.
Using that $u_{t_1}\nad v_{b+1}$ we can conclude that   $v_b\nad u_{t_1}$ as otherwise  $v_bv_{b+1}v_1v_{a-1}u_{t_1}v_b$ would be an induced 5-cycle in $G$.
  \qed

 In the case when $u_{q-1}$ is a cutvertex and hence $u_{t_k}=u_{q-1}$ we also
  cannot have $u_{t_1}\sim u_{t_k}$ as this would imply that $G$ contained the induced 5-cycle $u_{t_k}v_{b+1}v_1v_{a-1}u_{t_1}u_{t_k}$.

\begin{claim}
  \label{u2cl}
  $v_2\nad v_{b+1}$
\end{claim}
\pf
  Suppose that $v_2\sim v_{b+1}$.  Consider the two paths
  $R_1=v_{b+1}v_1v_{a-1}u_{t_1}$ and $R_2=v_{b+1}v_2u_{t_1}$ of lengths 3 and 2 respectively.  Let $P$ be any induced
  $(u_{t_1},u_{s_k})$-path in $G[u_{t_1},v_{b+1}]$. By Claim \ref{ut1nonadj}  $Pu_{t_k}v_{b+1}$ forms an induced cycle of length at least 5 in $G$ with either $R_1$ or $R_2$, contradiction. Note that $P$ may have length zero (if $k=2$ and $u_{t_1}=u_{s_k}$).
\qed

\begin{claim}
  \label{uskb+1}
   $u_{s_k}\sim v_{b+1}$
\end{claim}
\pf
  Suppose that  $u_{s_k}\nad v_{b+1}$.  Let $R_3=v_bv_{b+1}v_1v_{a-1}u_{t_1}$ and $R_4=v_{b}v_nv_2u_{t_1}$. Let $P'$ be any induced $(u_{t_1},u_{s_k})$-path in $G[u_{t_1},u_{s_k}]$. If $u_{s_k}\sim v_b$ then $P'v_b$ induce an odd cycle of length at least 5 with either $R_3$ or $R_4$. So assume that $u_{s_k}\nad v_b$. Then we are in the case when $u_{q-1}$ is a cutvertex and we have $u_{t_k}=u_{q-1}\sim v_{b+1}$. 
  Now we see that if $P'$ has even length then $P'u_{t_k}v_b$ induces and odd cycle of length at least 5 with $R_4$, contradiction. So $P'$ must have odd length. The predecessor $u^-_{q-1}$ of $u_{q-1}$ in the round ordering of $V(G)$ is a twin of $u_{q-1}$ in $H[2,q-1]$ and  $u^-_{q-1}\sim v_b$  as $u_{q-1}$ is not a pseudo-cutvertex of $G$. Furthermore,  as $u^-_{q-1}$ and $u_{q-1}$ are not true twins in $G$ there is a vertex $z$ in $G[b+1,n]$ such that $z\sim u_{q-1}$ and $z\nad u^-_{q-1}$. Now we see that $v_1v_{a-1}P'u^-_{q-1}v_bzv_1$ is an induced odd cycle,
  contradiction. Thus we must have  $u_{s_k}\sim v_{b+1}$.
\qed

\begin{claim}
  \label{vnnotsk}
  $v_n\nad u_{s_k}$
\end{claim}
\pf
  Suppose $v_n\sim u_{s_k}$. Then $u_{q-1}$ is not a cutvertex of $H$ so $u_{t_k}=v_b$. As $u_{t_1}\nad v_b$ by Claim \ref{ut1nonadj} we cannot have $u_{s_k}=u_{t_1}$. Let $R_5=u_{s_k}v_{b+1}v_1v_{a-1}u_{t_1}$ and $R_6=u_{s_k}v_nv_2u_{t_1}$ and let $P''$ be any induced $(u_{t_1},u_{s_k})$-path in $G[u_{t_1},u_{s_k}]$ then $P''$ induces an odd cycle of length at least 5 with one of $R_5,R_6$, contradiction.
\qed

\begin{claim}
  $v_{a-1}\nad v_n$.
\end{claim}
\pf
  Suppose $v_{a-1}\sim v_n$ and recall that we also have $v_n\sim v_b$ by (\ref{vav2vbvn}). By Claims \ref{ut1nonadj} and \ref{uskb+1} we also have $u_{t_1}\nad v_{b}, u_{t_1}\nad v_{b+1}$ and $u_{s_k}\sim v_{b+1}$. %This implies that we have $u_{s_k}\neq u_{t_1}$.
  As $v_2$ and $v_{a-1}$ are not true twins of $G$ there is a vertex $w$ in $G[v_{b+2},v_{n-1}]$ such that $v_{a-1}\nad w\sim v_2$. Now we see that the vertices $v_b,w,v_n,v_1,v_2,v_{a-1},u_{t_1}$ induce a copy of $F_3$, contradiction. So $v_{a-1}\nad v_n$.
\qed

Combining the claims above we see that $u_{s_k},v_b,v_{b+1},v_n,v_1,v_2,v_{a-1}$ induce a copy of $F_2$. This contradiction completes the proof when $u_2$ is a cutvertex.\\

\noindent{}{\bf Case 2.2: $v_2$ is not a cutvertex of $H$ but $v_{q-1}$ is a cutvertex.}

Denote $H''' = H[1,q-1]$. %As $v_1\nad u_{q-1}$ by (\ref{v1nad}) we see that $cc(H[1,q-1])\geq 2$. Then $H'''$ is one of the 4 types.
Let $H'''[s_i,t_i]$ ($1 \leq i \leq k$) be the canonical clique covering of $H'$. 
Then $u_{s_1} = v_1$ and $u_{t_k} = u_{q-1}$.
As $u_{q-1}$ is not a cutvertex of $G$, its predecessor $u^-_{q-1}$ in the order of $V(G)$ is adjacent to $u_q=v_b$ and it is a twin of $u_{q-1}$ in $G[1,b]$.
  As $u^-_{q-1}$ and $u_{q-1}$ are not true twins of $G$, there is a vertex  $v$ in $G[b+1,n]$ such that $u^-_{q-1}\nad v\sim u_{q-1}$. In particular $u^-_{q-1}\nad v_n$.
  \\

Suppose first that $k=2$. Then $H'''$ is of type 1 or type 4.

\begin{itemize}
\item When  $H'''$ is type 1 we have $u_{2}=v_{a-1}\sim u_{t_2}=u_{q-1}$ and then 
$v_nv_2v_au^-_{q-1}v_bv_n$ is an induced 5-cycle, contradiction. Note that we cannot have $v_2\sim u^-_{q-1}$  as that would imply that $v_2\sim u_{q-1}$, contradicting $\ell(u_{t_2})=v_{a-1}$, because $u^-_{q-1}$ is a true twin of $u_{q-1}$ in $G[1,b]$.
  \item Hence $H'''$ must be of type 4 and thus we have  $u_{s_2}=v_a$, $u_{t_2}=u_{q-1}$.
  \item We must have $v_{a-1}\sim v_{b+1}$ as $v_{b+1}v_1v_{a-1}v_au_{q-1}v_{b+1}$ is not an induced 5-cycle.
  \item %As $u_{q-1}$ and $u^ -_{q-1}$ are not true twins of $G$ there exists a vertex $v$ of $G[b+1,n]$ such that $u^-_{q-1}\nad v\sim u_{q-1}$.
    Now we see that  $vv_{a-1}v_au^-_{q-1}v_bv$ is an induced 5-cycle, contradiction ($v$ was defined above).
    
\end{itemize}

Assume now that $k\geq 3$. Consider first the case when $H'''$ is of one of the types 1,3,4. By Lemma \ref{lem:pathsintypes} there is an induced even $(u_{t_1},u_{t_k})$-path $P_1$ in $H'''$. \\

\begin{itemize}
\item We must have $v_{a-1}\sim v_{b+1}$ as otherwise $v_1P_1v_{b+1}v_1$ would be an induced odd cycle of length $\geq 5$ in $G$.
%\item As $v_1$ and $v_n$ are not true twins of $G$ we get that $v_n\sim v_b$.
\item If $H'''$ is type 1 or 4, then $v_nP_1[u_{t_1},u_{s_k}]u^-_{q-1}v_bv_n$ is an induced odd cycle of length $\geq 5$, so $H'''$ is of type 3.
%\item Then $v_n\sim u_{q-1}$ as otherwise $v_nP_1v_bv_n$ would be an induced odd cycle of length $\geq 5$.
\item It follows from Lemma \ref{lem:pathsintypes} that there is an induced $(u_{t_1},u_{s_k})$-path $Q$ of even length in $H'''$.
  \item Note that $u_{s_k}\nad v_b$ as $u_{q-1}$ is a cutvertex of $H$. Hence $u_{s_k}\nad v_{b+1}$ and therefore $v_{b+1}Qu_{q-1}v_{b+1}$ is an induced odd cycle of length $\geq 5$, contradiction.
\end{itemize}

The final case is when $H'''$ is type 2. Then $u_{s_2}=v_{a-1}$ so $u_{s_2+1}=v_a$.\\

\begin{itemize}
  \item By Lemma \ref{lem:pathsintypes} there exists an induced odd 
$(u_{s_2+1},u_{t_k})$-path $P_2$ in $H'''$.
\item We must have $v_{a-1}\sim v_{b+1}$ as otherwise $v_1v_{a-1}P_2v_{b+1}v_1$ would be an induced odd cycle of length $\geq 5$ in $G$.
%\item As usual this implies that $v_n\sim v_b$.
\item By Lemma \ref{lem:pathsintypes} there exists an induced $(u_{s_1},u_{t_k})$-path $P_3$ of even length.
\item We have $v_n\nad u_{t_k}$ since otherwise $v_nP_3[u_{t_1},u_{t_k}]v_n$ would be an odd induced cycle of length at least 5.
\item If  $u^-_{q-1}\sim v_{b+1}$ then $v_{b+1}P_3[u_{t_1},u_{s_k}]u^-_{q-1}v_{b+1}$ is an induced odd cycle of length at least 5 so $u^-_{q-1}\nad v_{b+1}$.
  \item Now $v_{b+1}v_2v_aP_3[u_{t_2},u_{s_k}]u^-_{q-1}v_bv_{b+1}$ is an induced odd cycle of length at least 5, contradiction.
%\item This implies that $v_n\sim u_{q-1}=u_{t_k}$ as otherwise $v_nP_3v_bv_n$ would be an induced odd cycle of length $\geq 5$.
 % \item Now $v_nP_3[u_{s_1},u_{s_k}]u^-_{q-1}u_{q-1}v_n$ is an oduced odd cycle of length $\geq 5$, contradiction.
  \end{itemize}

This final contradiction completes the proof of Theorem \ref{lem:nopseudocv}.\qed\\

Summarizing, we have the following:

\begin{tm} 
Let $G$ be a reduced proper circular-arc graph that is not a proper interval graph.
Then $G$ is locally complete 2-edge-colourable if and only if following hold:
\begin{itemize}
\item[(1)] $G$ is perfect.
\item[(2)] $G$ does not contain any of $F_1, F_2, F_3$ as an induced subgraph.
\item[(3)] $G$ has an even number of pseudo-cutvertices.
\item[(4)] Each weak block is either complete or of type 4.
\end{itemize}
\end{tm}
\pf This follows immediately from Theorems \ref{lem:weakblocks},
\ref{lem:onecutv} and \ref{lem:nopseudocv}.
\qed

%\section{Remarks and open problems}

%Since locally complete 2-edge-colourable graphs are claw-free, a locally complete
%2-edge-colourable graph whose complement is also locally complete 2-edge-colourable
%must be $\{K_{1,3},\overline{K_{1,3}}\}$-free. Using a result on
%$\{K_{1,3},\overline{K_{1,3}}\}$-free graphs in \cite{pst}, one can give 
%a description of all locally complete 2-edge-colourable graphs whose complements 
%are also locally complete 2-edge-colourable.

%\bibliography{refs}

\end{document}